\documentclass[10pt]{article}
\usepackage{amsmath,amssymb,amsfonts}
\usepackage[matrix,arrow,curve]{xy}
\author{Victor Tourtchine
\thanks{Partially supported by the
grants NSH-1972.2003.01, MK-451.2003.01}}
\title{On the other side of the bialgebra of chord diagrams}
\date{}

\newtheorem{theorem}{Theorem}[section]
\newtheorem{lemma}[theorem]{Lemma}
\newtheorem{statement}[theorem]{Statement}
\newtheorem{remark}[theorem]{Remark}
\newtheorem{definition}[theorem]{Definition}
\newtheorem{corollary}[theorem]{Corollary}
\newtheorem{proposition}[theorem]{Proposition}
\newtheorem{example}[theorem]{Example}

\def\mod{\mathop{\rm mod}\nolimits}
\newcommand\R{{\mathbb R}}
\newcommand\C{{\mathbb C}}

\newcommand\Z{{\mathbb Z}}
\newcommand\Q{{\mathbb Q}}
\newcommand\N{{\mathbb N}}
\newcommand\K{{\mathcal K}}
\newcommand\HH{{\mathcal H}}
\newcommand\ZZ{{\mathcal Z}}
\newcommand\kk{\Bbbk}

\newcommand\Endom{{\mathcal {END}}}
\newcommand\Assoc{{\mathcal {ASSOC}}}

\newcommand\BV{{\mathcal {BV}}}
\newcommand\Lie{{\mathcal {LIE}}}
\newcommand\Poiss{{\mathcal {POISS}}}
\newcommand\Gerst{{\mathcal {GERST}}}
\renewcommand\O{{\mathcal O}}

\newcommand\ie{{\it i.e. }}
\newcommand\cf{{\it cf. }}
\newcommand\etc{{\it etc}}

\newcommand\TD{\unitlength=0.2em
\begin{picture}(10,3)
 \put(0,0.1){\line(1,0){10}}
\end{picture}}

\newcommand\THREE{\unitlength=0.2em
\begin{picture}(40,9)
 \put(0,0){\line(1,0){40}}
 \qbezier(5,0)(15,9)(25,0)
 \qbezier(15,0)(25,9)(35,0)
 \qbezier(5,0)(20,15)(35,0)
\end{picture}}

\def\ad{{\rm ad}}

\def\nbox{\quad$\Box$}
\def\nboxm{\quad\Box}

\newcommand\rth{\refstepcounter{equation}}
\newcommand\numb{\rth{\rm \theequation}}
\numberwithin{equation}{section}

\begin{document}
\maketitle \sloppy

{\footnotesize
\begin{abstract}
In this paper we describe complexes whose homologies are naturally isomorphic to the first term  of the 
Vassiliev spectral sequence computing (co)homology of the spaces of long knots in $\R^d$, $d\ge3$. The first term of the Vassiliev
spectral sequence is concentrated in some angle of the second quadrant. In homological case the lower line of this term  is the bialgebra of chord diagrams
(or its superanalog if $d$ is even). We prove in this paper that the groups of the upper line are all trivial. In the same bigradings we compute the homology
groups of the complex spanned only by strata of immersions in the discriminant (maps having only self-intersections). We interprete the obtained groups as 
subgroups of the (co)homology groups of the double loop space of a $(d-1)$-dimensional sphere. In homological case the last complex is the  normalized
Hochschild complex of the  Poisson or Gerstenhaber (depending on parity of $d$) algebras operad. The upper line bigradings are spanned by the operad
of Lie algebras. To describe the cycles in these bigradings we introduce new homological operations on Hochschild complexes. These new
 operations are in fact the Dyer-Lashof operations induced by the action of the singular chains operad of little squares on Hochschild complexes.

\end{abstract}


\medskip

}
\setcounter{section}{-1}
\section{Introduction}\label{introduction}
\subsection{Spaces of long knots. Approach of V.~Vassiliev}
The {\it space of long knots} in $\R^d$ is the space of smooth embeddings $\R^1\hookrightarrow\R^d$ that
coincide with a fixed linear map $\R^1\hookrightarrow\R^d$ outside some compact set (depending on a knot). The
long knots form an open everywhere dense subset in the affine space $\K\simeq\R^\omega$ of all smooth maps
$\R^1\to\R^d$ with the same behavior at infinity. The complement $\Sigma$ of this dense subset is called the {\it
discriminant space}. It consists of the maps having self-intersections and/or singularities. Any cohomology class
$\gamma\in H^i(\K\setminus\Sigma)$ of the knot space can be realized as the linking coefficient with an
appropriate chain in $\Sigma$ of codimension $i+1$ in $\K$. In other words one has the Alexander duality:
$$
\tilde H^i({\K}\backslash\Sigma)\simeq \tilde H_{\omega -i-1}(\bar\Sigma),
\eqno(\numb)\label{eq01}
$$
where $\bar\Sigma$ designates the one-point compactification of the discriminant $\Sigma$.

Strictly speaking the isomorphism~\eqref{eq01} has no sense since $\omega$ is infinity. To define rigorously the
right-hand side of~\eqref{eq01} one needs to use finite-dimensional approximations of the space of long knots,
{\it i.~e.}  finite-dimensional spaces $\R^N\subset\K$ that are in general position with the discriminant
$\Sigma$,   \cf~\cite{V1} .

\subsection{Vassiliev's and Sinha's spectral sequences. Main results}
The main tool of Vassiliev's approach to computation of the (co)homology of the knot space is  {\it
simplicial resolution} $\sigma$ of the discriminant, {\it cf.}~\cite{V1}. The space $\sigma$ has a natural filtration
$$
\emptyset=\sigma_0\subset\sigma_1\subset\sigma_2\subset\dots.
\eqno(\numb)\label{eq02}
$$
Vassiliev conjectures that this filtration homotopically splits, \ie
$$
\bar\sigma\simeq\bigvee_{i=1}^{+\infty}(\bar\sigma_i/\bar\sigma_{i-1}).
\eqno(\numb)\label{eq03}
$$
This would imply the isomorphism
$$
\tilde H_*(\bar\Sigma)\equiv\tilde H_*(\bar\sigma)\simeq\bigoplus_{i=1}^{+\infty}\tilde
H_*(\bar\sigma_i/\bar\sigma_{i-1}), \eqno(\numb)\label{eq04}
$$
or in other words that the spectral sequence (called Vassiliev's main spectral sequence) associated with the
filtration~\eqref{eq02} stabilizes in the first term.

Filtration~\eqref{eq02} induces an increasing filtration in the homology groups of $\Sigma$, \ie  in the cohomology
of the knot space:
$$
H_{(0)}^*({\mathcal K}\backslash\Sigma)\subset H_{(1)}^*({\mathcal K}\backslash\Sigma)\subset H_{(2)}^*({\mathcal
K}\backslash\Sigma)\subset\dots,
\eqno(\numb\label{eq05})
$$
and decreasing dual filtration in the homology groups:
$$
H_*^{(0)}({\mathcal K}\backslash\Sigma)\supset H_*^{(1)}({\mathcal K}\backslash\Sigma)\supset H_*^{(2)}({\mathcal
K}\backslash\Sigma)\supset\dots.
\eqno(\numb\label{eq06})
$$

For $d\ge 4$ filtrations~\eqref{eq05},~\eqref{eq06} are finite for any dimension $*$. Vassiliev's main spectral
sequence in this case computes the graded quotient associated with these filtrations.

In the most intriguing case $d=3$ almost nothing is clear. However we can say something about the dimension $*=0$.
The filtration~\eqref{eq05} does not exhaust the whole cohomology of degree zero. The knot invariants obtained by
this method are called the Vassiliev invariants, or invariants of finite type. The dual space to the graded
quotient of the space of finite type knot invariants is the {\it bialgebra of chord diagrams}. The invariants and
the bialgebra in question were intensively studied in the last decade, \cf~\cite{BN}. The completeness conjecture for the
Vassiliev knot invariants is the question about the convergence of the filtration~\eqref{eq06} to zero for $d=3$,
$*=0$. The realization theorem of M.~Kontsevich~\cite{K} proves that the Vassiliev spectral sequence over $\Q$ for
$d=3$, $*=0$ computes the corresponding associated quotient (for positive dimensions $*$ in the case $d=3$
even this is not for sure) and does stabilize in the first term. The groups of the graded quotient associated to
filtration~(\ref{eq05}) in the case $d=3$, $*>0$ are some quotient groups of the groups calculated by Vassiliev's
main spectral sequence.

\begin{figure}[!ht]
\begin{center}
\unitlength=0.5em
\begin{picture}(20,30)
\put(0,10){\vector(1,0){20}}
\put(10,0){\vector(0,1){30}}
\put(0,30){\line(1,-2){10}}
\put(5,30){\line(1,-4){5}}
\put(19,8.5){$p$}
\put(10.5,29){$q$}
\put(8,14){\line(1,0){1}}
\put(6,18){\line(1,0){2}}
\put(4,22){\line(1,0){3}}
\put(2,26){\line(1,0){4}}
\end{picture}
\end{center}
\caption{First term of the Vassiliev's spectral sequence}\label{SpSeq}
\end{figure}
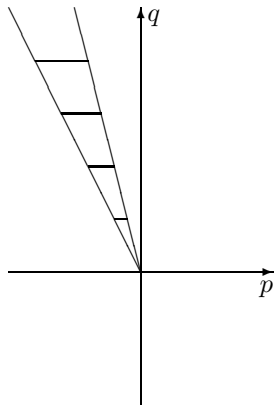

The first term of  Vassiliev's spectral sequence is concentrated in some angle of the second quadrant, see Figure~\ref{SpSeq}.

The lower and the upper edges of this angle are determined by the following equations:
$$
\begin{tabular}{ll}
$q=-(d-2)p$&\text{lower line;}\\
$q=-(d-1)p-1$&\text{upper line.}
\end{tabular}\eqno(\numb\label{eq07})
$$
In the cell $(p,q)$ stands the group $\tilde H^{\omega-p-q-1}(\bar\sigma_{-p}/\bar\sigma_{-p-1})$. Observe, that
$p$ is equal to {\it minus} the filtration: $p=-i$, and $q$ is defined by the condition that $p+q$ equals  the usual
(co)homology dimension.

Up to a changing of grading this first term depends only on the parity of $d$ (dimension of the ambient space
$\R^d$). For different $d_1$ and $d_2$ of the same parity one has the isomorphism:
\begin{align*}
_{d_1}E_1^{p,q}&\simeq  {\, _{d_2}E_1^{p,q-p(d_2-d_1)}},\\
_{d_1}E^1_{p,q}&\simeq  {\, _{d_2}E^1_{p,q-p(d_2-d_1)}},
\end{align*}
where $_dE_1^{p,q}$ (cohomological case), $_dE^1_{p,q}$ (homological case) designate the first term of the Vassiliev spectral 
sequence computing the (co)homology of  the space of long knots  in $\R^d$. 
Over $\Z_2$ this is also true for  dimensions $d_{1}$, $d_{2}$ of different parities.  In the homological case the direct sum of groups
standing in the lower line is isomorphic to the bialgebra of chord diagrams, if $d$ is odd, and to a non-trivial superanalog of
this bialgebra, if $d$ is even.

Actually there is another and absolutely different approach to studing the space of embeddings. Briefly speaking in this approach one  
"approximates"  the space of knots $Emb=\K\setminus\Sigma$ by means of homotopy limits of diagrams of maps. This approach was initiated by 
T.~Goodwillie and M.~Weiss~\cite{Good,GW}, and then developped by D.~Sinha~\cite{Sinha1,Sinha2}, and also by I.~Volic~\cite{Vol}.
 In particular for the space of long knots this method provides a spectral sequence whose second term is isomorphic up to a shift of bigradings to
 the first term  of the Vassiliev spectral sequence (this spectral sequence was constructed by D.~Sinha~\cite{Sinha1}):  
 
\begin{proposition}\label{p0}
{\rm (i)} The groups of the first term of Vassiliev's spectral sequence (computing the (co)homology of the space of long knots) are naturally 
isomorphic up to a shift of bigradings to the groups of the second term of Sinha's spectral sequence 
(computing the (co)homology of a space homotopy equivalent to the space of long knots).

{\rm (ii)}  Moreover, the first term of Vassiliev's auxiliary spectral (stabilizing in the second term and computing the first term of the main spectral sequence) 
is isomorphic as a complex to the first term of Sinha's spectral. \nbox
\end{proposition}
 
Assertion (ii) says that  the complexes $CT_{*}D^{odd}$, $CT_{*}D^{even}$ that
 we define in Section~\ref{s3} is nothing else but the first term of Sinha's spectral sequence. To prove Proposition~\ref{p0} one  needs to compare
 Section~7 in~\cite{Sinha1} with Vassiliev's auxiliary spectral sequence which was  defined in the seminal work~\cite{V1}. For a more explicit description of
 the auxiliary spectral sequence, see also~\cite{V4} or~\cite[Chapitre~II]{T2}.

It is also worth to mention the Cattaneo-Cotta-Ramusino-Longoni  construction~\cite{CCL,CCL2} that 
provides a morphism from a graph-complex to the De Rham complex
of the space of knots. In the case of long knots the corresponding graph-complex is quasi-isomorphic to the complex $CT_{*}D$.

\medskip

In our paper we find the groups of the upper line~\eqref{eq07}. Actually we prove that all these groups are
trivial.

We also consider  complexes that arise if in the degree zero term  we take into account only the summands
corresponding to strata assigned to maps having only self-intersections (but not degeneration of the derivative).
For these complexes we find the homology groups in the bigradings of the upper line. 
I.~Volic and D.~Sinha give a geometrical interpretation of these complexes. In their approach the complexes are the first terms of the 
spectral sequences computing the cohomology groups of the homotopy 
fiber of the inclusion of the space of knots $Emb=\K\setminus\Sigma$ to the space of immersions $Imm$, \cf~\cite{Vol,Sinha2}.  
D.~Sinha proved that this homotopy fiber is hopotopy equivalent to  a direct product $Emb\times \Omega^2 S^{d-1}$, see~\cite{Sinha2}.
Therefore, the non-trivial upper diagonal homology groups that we find in this paper are some subgroups of the (co)homology of $\Omega^2 S^{d-1}$.
One more reason why the above
complexes are worth being studied is that in homological case they are isomorphic to the normalized Hochschild
complexes of the Poisson algebras operad ($d$ odd) and Gerstenhaber algebras operad ($d$ even), \cf~\cite{T1,T2,T4}. Their
homology groups are the characteristic classes of Hochschild cohomology of Poisson, resp. Gerstenhaber algebras
considered as associative algebras. The homology groups of the upper line bigradings are the characteristic  classes
defined by the Lie-algebra structure in Poisson or Gerstenhaber algebras.

\subsection{Plan of the paper}
The paper is divided  into 4 parts. Each part starts with a brief description of its contents. The results of Parts~\ref{DHATD} and~\ref{s8} 
were given before in the thesis~\cite{T2}. The main results of Part~\ref{s8} were given also (and with more details) in~\cite{T4}.
The results of Part~\ref{upper_diagonal} and of Appendix~\ref{E} are new.

Part~\ref{DHATD} is cohomological: in this part we define comlexes computing the first term of the Vassiliev spectral sequence converging 
to the cohomology of the space of long knots. We remind that one conjectures that this spectral sequence stabilizes in the first term.

Part~\ref{s8} is homological: in this part  we define complexes that are dual to those of Part~\ref{DHATD}. We give a hint how these
dual complexes were obtained and describe their bases. The obtained complexes are the normalized Hochschild complexes of some
operads endowed with a map from the associative algebras operad.

In Part~\ref{upper_diagonal} we study the upper diagonal homology groups of the complexes defined in Parts~\ref{DHATD}-\ref{s8}. In this part 
we define new homological operations on Hochschild complexes. It will be proven elsewhere that these operations are the Dyer-Lashof homology
operations induced by the little square chains action on Hochschild complexes, see also Appendix~\ref{E}. The new results in this part 
are given by Theorems~\ref{t71}, \ref{t72}, \ref{t116}. Theorem~\ref{t116} is given without proof.

The last part contains various appendixes. It is a mixture of explicit demonstration of objects we deal with; results of computer calculations; 
and reformulation of old results with some nuances that are necessary for our considerations. In Appendix~\ref{E} 
we show how the homology bialgebra of $\Omega^2S^{d-1}$ is included in the Hochschild homology of Poisson and Gerstenhaber 
algebras operads. This result is new: we use there our new "Dyer-Lashof" operations on Hochschild complexes; but it was more logical to 
put this section after Appendix~\ref{D}.

\medskip

All along the paper we make a confusion of $\Z$-grading $deg=p+q$  and corresponding $\Z_{2}$ supergrading. It is done delibirately
to emphasize the fact that the obtained complexes depends only on the parity of  dimension $d$ of the ambient space $\R^d$.

\subsection{Acknowledgement} 
I whould like to thank Universit\'e Catholique de Louvain (Louvain-la-Neuve) where this paper was finished for hospitality.

I am grateful to V.~Vassiliev, P.~Lambrechts, D.~Chataur and G.~Sharygin for interesting conversations. 
I am deeply grateful to J.-O.~Moussafir and A.~Semionov for explaining me many things in C++. 
I thank J.-L.Loday, J.~McClure, D.~Sinha and A.~Voronov for giving me different
references and encouragement. I am also indebted to G.~Sharygin for correcting the earlier version of this paper.

\part{Differential Hopf algebras of $T_*$-diagrams, $T$-diagrams and $T_0$-diagrams}\label{DHATD}
To compute the first term of the main spectral sequence V.~A.~Vassiliev introduced an {\it auxiliary filtration}
on the spaces $\sigma_i\setminus\sigma_{i-1}$. The associated {\it auxiliary} spectral sequence stabilizes in the
second term since its first term is concentrated on the only row. The degree zero term with the differential on it
is a direct sum of tensor products of so-called {\it complexes of connected graphs}. The homology groups of the complex
of connected graphs on any finite set $M$ are trivial everywhere except the minimal possible dimension. This only non-trivial group
is the $\Z$-module $T_M^-$ described in Section~\ref{s21}. Complex $CT_*D$ of $T_*$-diagrams that we define in
Section~\ref{s3} is exactly the first term of the auxiliary spectral sequence in cohomological case --- \lq\lq
cohomological\rq\rq means relation to the \underline{cohomology} of the knot space. 

In Section~\ref{s4} we define a base in the space of the complex of $T_*$-diagrams.

In Section~\ref{s5} we introduce the complexes of $T$-diagrams and of $T_0$-diagrams. The first complex $CTD$ is a
quotient-complex of  complex $CT_*D$. It is spanned by strata of maps without degeneration of derivative. The
complex $CT_0D$ of $T_0$-diagrams is a subcomplex of $CTD$ and $CT_*D$. It is quasi-isomorphic to the complex
$CT_*D$ and serves to simplify computations of the homology groups of $CT_*D$:
$$
\xymatrix{
CT_{0}D\ar@{^{(}->}[rr]^{\sim }\ar@{^{(}->}[dr]&&CT_{*}D\ar@{->>}[dl]\\
&CTD&
}
\eqno(\numb\label{eq06'})
$$

In Section~\ref{s5} we describe the differential Hopf algebra structure on the complexes $CT_*D$, $CTD$, $CT_0D$.
The corresponding differential Hopf algebras are designated by $DHAT_*D$, $DHATD$, $DHAT_0D$. So, \eqref{eq06'} is a commutative 
diagram of differential Hopf algebras morphisms.

\section{$\Z$-modules $T_M^-$, $T_M^+$. Complexes of $T_*$-diagrams}\label{s2}
\subsection{$\Z$-modules $T_M^-$, $T_M^+$}\label{s21}
\subsubsection{$T_M^-$}
Consider a finite set $M$ of some cardinality $\#M$.  We will define an {\it orientation} of a tree with $\#M$
vertices labelled  one-to-one by the elements of $M$ as an ordering of its edges. Consider a $\Z$-module spanned
by the oriented trees --- changing of order of the edges implies multiplication of such a tree by
$(-1)^{|\sigma|}$, where $|\sigma|$ is the parity of the corresponding permutation. An oriented tree can be
viewed as a monomial of anticommuting elements $\alpha_{ab}=\alpha_{ba}$ representing edges ($a,b\in M$ and
$a\ne b$).

$\Z$-module $T_M^-$ is defined as the quotient-space of the above $\Z$-module by all the 3-term relations of the
following type:

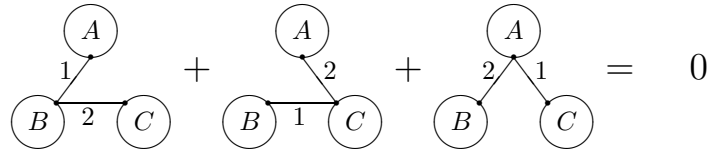
\begin{figure}[!ht]
\unitlength=0.2em
\begin{center}
\begin{picture}(126,33)
 \put(0,0){
   \begin{picture}(30,27.16)
      \put(5,5){\circle{10}}
       \put(8.52,8.52){\circle*{1.2}}
      \put(25,5){\circle{10}}
      \put(21.48,8.52){\circle*{1.2}}
      \put(15,22.16){\circle{10}}
      \put(15,17.16){\circle*{1.2}}
      \put(8.52,8.52){\line(3,4){6.48}}
      \put(8.52,8.52){\line(1,0){12.96}}
      \put(3.1,3.5){$B$}
      \put(23.1,3.5){$C$}
      \put(13.1,20.66){$A$}
      \put(9,13){$1$}
      \put(13.2,4.3){$2$}
      \end{picture} 
    }
   \put(34,13){\Large $+$}   
   \put(40,0){
      \begin{picture}(30,27.16)
      \put(5,5){\circle{10}}
       \put(8.52,8.52){\circle*{1.2}}
      \put(25,5){\circle{10}}
      \put(21.48,8.52){\circle*{1.2}}
      \put(15,22.16){\circle{10}}
      \put(15,17.16){\circle*{1.2}}
      \put(15,17.16){\line(3,-4){6.48}}
      \put(8.52,8.52){\line(1,0){12.96}}
      \put(3.1,3.5){$B$}
      \put(23.1,3.5){$C$}
      \put(13.1,20.66){$A$}
      \put(19,13){$2$}
      \put(13.2,4.3){$1$}
      \end{picture} 
    }
    
     \put(74,13){\Large $+$}   
     \put(80,0){
   \begin{picture}(30,27.16)
      \put(5,5){\circle{10}}
       \put(8.52,8.52){\circle*{1.2}}
      \put(25,5){\circle{10}}
      \put(21.48,8.52){\circle*{1.2}}
      \put(15,22.16){\circle{10}}
      \put(15,17.16){\circle*{1.2}}
      \put(8.52,8.52){\line(3,4){6.48}}
      \put(15,17.16){\line(3,-4){6.48}}
      \put(3.1,3.5){$B$}
      \put(23.1,3.5){$C$}
      \put(13.1,20.66){$A$}
      \put(9,13){$2$}
      \put(19,13){$1$}           
      \end{picture} 
    }  
    \put(114,13){\Large $=$ \hspace{4mm} $0$}
\end{picture}
\end{center}
\caption{3-term relations in $T_M^-$}\label{fig1}
\end{figure}

This picture is a sum of 3 trees whose edges are the same except the first two. This relation in terms of monomials
 can be rewritten as
$$
(\alpha_{ab}\alpha_{bc}+\alpha_{bc}\alpha_{ca}+\alpha_{ca}\alpha_{ab})\cdot T_A\cdot T_B\cdot T_C=0,
\eqno(\numb\label{3T})
$$
where $M=A\sqcup B\sqcup C$, $a\in A$, $b\in B$, $c\in C$; $T_A$, $T_B$ and $T_C$ are some trees on sets
 $A$, $B$ and
$C$ respectively.

In this form the relations resemble the Arnold's relations in the cohomology algebra of the configuration spaces
of different points in $\R^2$, \cf~\cite{A1}. Indeed $T_{M}^-$ is the maximal degree cohomology group 
 of the space of configurations of $\# M$ points (labelled by the elements of $M$).

\begin{statement}\label{st1}
{\rm \cite{V3}} $T_M^-\simeq\Z^{(\#M-1)!}$. \nbox
\end{statement}

In the case $M=\{1,2,\dots,n\}$ the module $T_M^-$ will be also denoted by $T_n^-$. Note that $T_n^-$ is endowed
with an action of the symmetric group $S_n$.

\subsubsection{$T_M^+$}
For the same finite set $M$ we define another $\Z$-module $T_M^+$ that is also spanned by the trees and quotiented
by some 3-term relations. We will orientate trees in another way. By definition an orientation of a tree is an
orientation of all its edges. Note we demand no more their ordering. Changing of orientation of one of edges is
equivalent to multiplication by $(-1)$. An oriented tree can be viewed as a monomial of commuting elements
$\alpha_{ab}=-\alpha_{ba}$ representing oriented edges ($a,b\in M$ and $a\ne b$).

The 3-term relations are as follows:

\begin{figure}[!ht]
\unitlength=0.2em
\begin{center}
\begin{picture}(126,33)
 \put(0,0){
   \begin{picture}(30,27.16)
      \put(5,5){\circle{10}}
       \put(8.52,8.52){\circle*{1.2}}
      \put(25,5){\circle{10}}
      \put(21.48,8.52){\circle*{1.2}}
      \put(15,22.16){\circle{10}}
      \put(15,17.16){\circle*{1.2}}
      \put(15,17.16){\vector(-3,-4){6.4}}
      \put(8.52,8.52){\vector(1,0){12.8}}
      \put(3.1,3.5){$B$}
      \put(23.1,3.5){$C$}
      \put(13.1,20.66){$A$}
       \end{picture} 
    }
   \put(34,13){\Large $+$}   
   \put(40,0){
      \begin{picture}(30,27.16)
      \put(5,5){\circle{10}}
       \put(8.52,8.52){\circle*{1.2}}
      \put(25,5){\circle{10}}
      \put(21.48,8.52){\circle*{1.2}}
      \put(15,22.16){\circle{10}}
      \put(15,17.16){\circle*{1.2}}
      \put(8.52,8.52){\vector(1,0){12.8}}
      \put(21.48,8.52){\vector(-3,4){6.4}}
      \put(3.1,3.5){$B$}
      \put(23.1,3.5){$C$}
      \put(13.1,20.66){$A$}
  
      \end{picture} 
    }
    
     \put(74,13){\Large $+$}   
     \put(80,0){
   \begin{picture}(30,27.16)
      \put(5,5){\circle{10}}
       \put(8.52,8.52){\circle*{1.2}}
      \put(25,5){\circle{10}}
      \put(21.48,8.52){\circle*{1.2}}
      \put(15,22.16){\circle{10}}
      \put(15,17.16){\circle*{1.2}}
      \put(21.48,8.52){\vector(-3,4){6.4}}
      \put(15,17.16){\vector(-3,-4){6.4}} 
      \put(3.1,3.5){$B$}
      \put(23.1,3.5){$C$}
      \put(13.1,20.66){$A$}
      \end{picture} 
    }  
    \put(114,13){\Large $=$ \hspace{4mm} $0$}
\end{picture}
\end{center}
\caption{3-term relations in $T_M^+$}\label{fig2}
\end{figure}
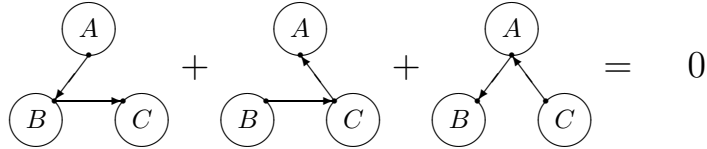

In terms of monomials this relation is given by the same formula~\eqref{3T}.

Let us denote by $\pm\Z$ the one-dimensional sign representation of the symmetric group $S_n$. Let $T^+_n$
designate $T_{\{1,2,\dots,n\}}^+$.

\begin{statement}\label{st2}
{\rm \cite{T2}}
There is a natural isomorphism of the $S_n$-modules $T_n^+\simeq T_n^-\otimes (\pm\Z).$ \nbox
\end{statement}

In particular this statement implies $T_M^+\simeq\Z^{(\#M-1)!}$.

\subsection{Complexes of $T_*$-diagrams}\label{s3}
Now we are ready to describe the complexes of $T_*$-diagrams. We denote them $CT_*D^{odd}$ in the case of odd $d$
and $CT_*D^{even}$ in the case of even $d$. When we want to treat them simultaneously we will write simply
$CT_*D$.  We remind that these complexes compute the first term of the Vassiliev spectral sequence, see
Figure~\ref{SpSeq}.
The complexes are bigraded. Bigrading $(i,j)$ corresponds to bigrading $(p,q)=(-i,id-j)$ of the Vassiliev spectral sequence.

\subsubsection{case of even $d$}
Consider a finite set of points on the line $\R^1$. We define a $T_*$-diagram on this set as a number of edges
joining its points and also a number of asterisks that are  put in these points.

We demand

\nopagebreak

1) each point contains not more then 1 asterisk;

2) the resulting graph is a disjoint union of trees;

3) if  a tree consists of only one point, this point must contain an asterisk.

\medskip

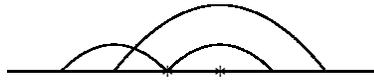
\begin{figure}[!ht]
\unitlength=0.2em
\begin{center}
\begin{picture}(70,15)
 \put(0,0){\line(1,0){70}}
 \qbezier(10,0)(20,10)(30,0)
 \qbezier(30,0)(40,10)(50,0)
 \qbezier(20,0)(40,25)(60,0)
 \put(29,-1.2){$*$}
 \put(39,-1.2){$*$}
\end{picture}
\end{center}
\caption{example of a $T_*$-diagram}\label{fig3}
\end{figure}

Figure~\ref{fig3} gives an example of a $T_*$-diagram.

The points of a $T_*$-diagrams are split into a partition. The subsets of this partition formed by the vertices of
the trees will be called {\it minimal components}. For example, the $T_*$-diagram on the Figure~\ref{fig3} has 3
minimal components: $\{1,3,5\}$, $\{2,6\}$ and $\{4\}$.

An {\it orientation} of a $T_*$-diagram is by definition an ordering of its {\it orienting set} that consists of
three types of elements:

\begin{center}
1) points; \quad 2) edges;\quad 3) asterisks.
\end{center}

Changing of order is equivalent to multiplication by the sign of the corresponding permutation.

\bigskip

For example, the orientation set of the diagram from Figure~\ref{fig3} has 11 elements:

$\alpha_{13}$, $\alpha_{35}$, $\alpha_{26}$ --- representatives of edges;

$\beta_1$, $\beta_2$, $\beta_3$, $\beta_4$, $\beta_5$, $\beta_6$   --- representatives of points;

$\gamma_{3}^* $,  $\gamma_{4}^* $ --- representatives of asterisks.

\medskip

To orient the diagram one needs to order these elements, for example as it is done below:

$$
\beta_1\beta_2\beta_3\beta_4\beta_5\beta_6\alpha_{13}\alpha_{35}\alpha_{26}
\gamma_{3}^*\gamma_{4}^* .
$$

We can consider this listing as a monomial of anticommuting elements.

\bigskip

Diagrams obtained one from another by an orientation preserving diffeomorphism  of the line are set to be
equivalent or equal.

The space of the complex of $T_*$-diagrams is defined as the $\Z$-module spanned by the oriented $T_*$-diagrams
and quotiented by all the 3-term relations issued from the relations in the $\Z$-modules $T_M^-$, $M$ being any
minimal component. This space is a direct sum of tensor products of modules $T_M^-$.

\bigskip
\bigskip

Now let us describe the differential in this complex $CT_*D^{even}$.

First of all note that the space of this complex is bigraded:

1) the first grading $i$ --- we call it {\it complexity} --- is the total number of edges and stars of a diagram;

2) the second grading $j$ is the number of points of a diagram.

\medskip

The differential $\partial$ will be of bigrading $(0,-1)$. It conserves the complexity but diminishes by one the 
number of
points on the line. It means that complex $CT_*D$ is a direct sum over $i$ of complexes, each of them being
spanned by the diagrams of complexity $i$.

The differential $\partial$ of a $T_*$-diagram $D$ with $j$ points is a sum over all possible $j-1$ gluings of two
neighbor points on the line $\R^1$. A gluing that does not give zero will be called {\it admissible}. There are
two possibilities:

1) the gluing points $t_1$ and $t_2$ belong to different minimal components;

2) the gluing points $t_1$ and $t_2$ belong to the same minimal component.

\medskip

In the first case the gluing is admissible if and only if there is no asterisk in at least one of these two
points. The boundary diagram is obtained from $D$ by contraction of the segment $[t_1,t_2]$, see
Figure~\ref{fig4}:

\begin{figure}[!ht]
\unitlength=0.2em
\begin{center}
\begin{picture}(120,40)
 \put(0,0){
  \begin{picture}(30,40)
   \put(0,25){
     \begin{picture}(30,15)
       \put(0,0){\line(1,0){30}}
       \qbezier(3,10)(10,8)(12,0)
       \qbezier(27,10)(20,8)(18,0)
     \end{picture}
     }
 
   \put(16.5,22){\vector(0,-1){10}}
   
   \put(0,0){
     \begin{picture}(30,15)
       \put(0,0){\line(1,0){30}}
       \qbezier(3,10)(10,8)(15,0)
       \qbezier(27,10)(20,8)(15,0)
     \end{picture}
     }
    \end{picture} 
  }
  
  \put(50,0){
     \begin{picture}(70,40)
        \put(0,25){
        \begin{picture}(30,15)
           \put(0,0){\line(1,0){30}}
           \qbezier(3,10)(10,8)(12,0)
           \qbezier(27,10)(20,8)(18,0)
           \put(10.8,-1.1){$*$}
        \end{picture}
        }
        
        \put(40,25){
        \begin{picture}(30,15)
           \put(0,0){\line(1,0){30}}
           \qbezier(3,10)(10,8)(12,0)
           \qbezier(27,10)(20,8)(18,0)
           \put(16.8,-1.1){$*$}
        \end{picture}
        }
        
        \put(20,0){
        \begin{picture}(30,15)
          \put(0,0){\line(1,0){30}}
          \qbezier(3,10)(10,8)(15,0)
          \qbezier(27,10)(20,8)(15,0)
           \put(13.9,-1.15){$*$}
        \end{picture}
        }
     
     \put(19,22){\vector(1,-1){7}}
     \put(54.5,22){\vector(-1,-1){7}}
     \end{picture}
 }    
  
\end{picture}
\end{center}
\caption{admissible gluings of points  from different minimal components}\label{fig4}
\end{figure}
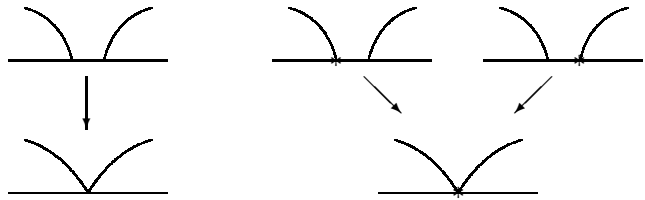

To obtain the orienting monomial (with some sign) of this boundary diagram one  needs to place on the first 
place of the orienting set of $D$ the representative of the left point $t_1$, and on the second place ---
representative of  the right point $t_2$ (this gives a sign) and then replace them by one element --- representative 
of the new point of contraction. All the rest in the orienting monomial do not change.

\medskip

In the second case ($t_1$ and $t_2$ from the same minimal component) gluing is admissible if and only if the
points $t_1$ and $t_2$ are joined by an edge and none of them contains an asterisk. The boundary diagram is
obtained by contracting  the segment $[t_1,t_2]$ and replacing the edge $\alpha_{t_{1}t_{2}}$ joining $t_1$ 
with $t_2$ by new asterisk $\gamma^*_{t}$, where $t$ is the point of contraction,
see Figure~\ref{fig5}:

\begin{figure}[!ht]
\begin{center}
  \unitlength=0.2em
  \begin{picture}(80,15)
    \put(50,0){
        \begin{picture}(30,15)
          \put(0,0){\line(1,0){30}}
          \qbezier(3,10)(10,8)(15,0)
          \qbezier(27,10)(20,8)(15,0)
           \put(13.9,-1.15){$*$}
        \end{picture}
        }
        
     \put(36,6){\vector(1,0){10}}   
     \put(0,0){
       \begin{picture}(30,15)
        \put(0,0){\line(1,0){30}}
        \qbezier(3,10)(10,8)(12,0)
        \qbezier(27,10)(20,8)(18,0)
        \qbezier(12,0)(15,5)(18,0)
     \end{picture}
     }
    
  \end{picture}
\end{center}
\caption{admissible gluing of points from the same minimal component}\label{fig5}
\end{figure}
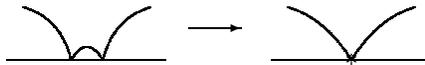

The orientation monomial (with some sign) of the boundary diagram is obtained analogously.

\subsubsection{case of odd $d$}

Now we will define the complex $CT_*D^{odd}$ of $T_*$-diagram in the case of odd $d$. The space of this complex
is also defined to be spanned by $T_*$-diagrams and quotiented by 3-term relations. The difference is that we have
another definition of orientation of $T_*$-diagrams. The orienting set of a diagram contains now only its points on the line.
For example the orienting set of the diagram from  Figure~\ref{fig3} consists of 6 elements:
$\beta_1$, $\beta_2$, $\beta_3$, $\beta_4$, $\beta_5$, $\beta_6$   --- representatives of its 6 points; as its orienting 
monomial we can thus take
$$
\beta_1\beta_2\beta_3\beta_4\beta_5\beta_6.
$$
As in the case of odd $d$ we consider orienting elements as anticommuting generators.
To orient a $T_*$-diagram one should also fix orientation of all its edges. Changing of orientation of one of
edges is equivalent to multiplication by $(-1)$.
 For instance, one can orient all
edges from left to right similarly to the orientation of $\R^1$:

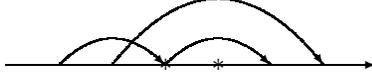
\begin{center}
\unitlength=0.2em
\begin{picture}(70,15)
 \put(0,0){\vector(1,0){70}}
 \qbezier(10,0)(20,10)(30,0)
 \qbezier(30,0)(40,10)(50,0)
 \qbezier(20,0)(40,25)(60,0)
 \put(29,-1.2){$*$}
 \put(39,-1.2){$*$}
 \put(30,0){\vector(1,-1){0}}
 \put(50,0){\vector(1,-1){0}}
 \put(60,0){\vector(2,-3){0}}
\end{picture}
\end{center}

3-term relations in this space are issued from those in the $\Z$-modules $T_M^+$, $M$ being any minimal component.
So, the space of complex $CT_*D^{odd}$ is a direct sum of tensor products of modules $T_M^+$.

The differential in this complex is defined analogously  to the case of even $d$:  as the sum of admissible gluings with appropriate 
signs. When we take a gluing of two points joined by an edge whose orientation is opposit to that of the line, we should additionally multiply 
the result of this gluing by minus~1.

\section{Basis in the complexes $CT_*D$. Alternated $T_*$-diagrams}\label{s4}
\subsection{Basis in $T_M^-$, $T_M^+$}
There exist three different types of basis for these $\Z$-modules: basis of snakes, \cf~\cite{V3}, basis of monotone
trees and basis of alternated trees, \cf~\cite{T2}.  To define the first basis one needs to fix one element in $M$ ---
the head of the snakes. To define the second and the third bases one needs the elements of $M$ to be linearly ordered.
Actually, only the third type of basis turned out to be  useful  in the investigation of complex $CT_*D$
and its derivatives $CTD$ and $CT_0D$, defined in the next section. The reason for this is that only this basis is preserved
via inversion of the order in $M$.

\smallskip
Let finite set $M$ be linearly ordered.
We will consider trees, whose vertices are labeled one-to-one by elements of $M$, as rooted trees with a root being the minimal element of $M$.

\begin{definition}\label{partialOrder}
{\rm Any rooted tree defines a {\it partial order} $\prec$ on the set $M$ of its vertices: we say that $a\prec b$ for two different elements of $M$ if and only if
the (only) path from the root to $b$ passes through $a$.} \nbox
\end{definition}

\begin{definition}\label{alt_trees}
{\rm For a linearly ordered set $(M,<)$ a tree, whose vertices are labelled one-to-one by the elements of $M$, is called {\it
alternated} if for any its edge $(a,b)$ and for any element $x\in M$ such that $a\prec x$, $b\prec x$, one always has $a<x<b$ or $b<x<a$.}
\nbox
\end{definition}

Note that all  alternated trees contain the edge joining the  extremal vertices.

One can give an equivalent recursive definition:

\begin{definition}\label{alt_trees1}
{\rm A tree, whose vertices are linearly ordered, is called {\it
alternated} if it is a trivial tree with only one vertex or it contains the edge joining the extremal vertices and when one removes this edge the remaining
two disconnected trees are alternated.}
\nbox
\end{definition}

\begin{figure}[!ht]
\begin{center}
\unitlength=0.2em
\begin{picture}(50,23)(0,-3)
 \qbezier(0,0)(15,15)(30,0)
 \qbezier(0,0)(25,30)(50,0)
 \qbezier(10,0)(30,20)(50,0)
 \qbezier(10,0)(25,15)(40,0)
 \qbezier(20,0)(30,10)(40,0)
 \multiput(0,0)(10,0){6}{\circle*{1.2}}
 \put(-1.5,-5){$1$}
 \put(8.5,-5){$2$}
 \put(18.5,-5){$3$}
 \put(28.5,-5){$4$}
 \put(38.5,-5){$5$}
 \put(48.5,-5){$6$}
\end{picture}
\end{center}
\caption{an alternated tree}\label{fig5'}
\end{figure}
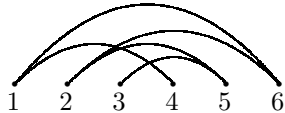

\begin{lemma}\label{l41}
{\rm \cite{T2}}
For any linearly ordered set $M$ of cardinality $\#M$ there are exactly $(\#M-1)!$ alternated trees on it, and they form a
basis in the $\Z$-modules $T_M^-$, $T_M^+$. \nbox
\end{lemma}

\subsection{Alternated $T_*$-diagrams}

\begin{definition}\label{alt_diagr}
{\rm A $T_*$-diagram is called {\it alternated} if all its trees are alternated, where the order on the minimal components
is induced by the order on the line $\R^1$.} \nbox
\end{definition}

\section{Complexes of $T$-diagrams and $T_0$-diagrams}\label{s5}

\begin{definition}\label{t-diagr}
{\rm A $T_*$-diagram is called {\it $T$-diagram} if it does not contain asterisks. \nbox}
\end{definition}

\begin{definition}\label{t0-diagr}
{\rm A $T$-diagram is called {\it $T_0$-diagram} if it does not have edges joining two neighbor points. \nbox}
\end{definition}

The subspace spanned by $T$-diagrams forms a quotient-complex of complex $CT_*D$. We will call it {\it complex
of $T$-diagrams} or simply $CTD$:
$$
\xymatrix{
CT_*D\ar@{->>}[r]&CTD.
}
$$
We quotient by the subspace spanned by the diagrams having asterisks. Alternated $T$-diagrams form a basis in this complex.

\begin{proposition}\label{p51}
The space spanned by $T_0$-diagrams forms a subcomplex of $CT_*D$ (and therefore of $CTD$). This inclusion is
a quasi-isomorphism of complexes. The set of alternated $T_0$-diagrams forms a basis in $CT_0D$. \nbox
\end{proposition}

\noindent {\bf Proof:}  It is not evident from Definition~\ref{t0-diagr} that the space of $T_{0}$-diagrams forms a subcomplex 
of $CT_{*}D$. It will be clear from the proof that this space is actually the maximal subspace of $CT_{*}D$
that lies in the space of $T$-diagrams and is invariant with respect to the differential.

Consider the spectral sequence associated to the filtration in $CT_*D$ by the number of
minimal components. The degree zero differential $d_0$ of this spectral sequence is the sum of admissible gluings
of the second type (gluing of points from the same minimal component), see Section~\ref{s3}. If we examine the differential
$d_0$ in the basis of alternated $T_*$-diagrams, we obtain that the first term of this spectral sequence is
concentrated on the only row and is spanned by the alternated $T_0$-diagrams. This proves the proposition. \nbox

\bigskip
\bigskip

One gets the following commutative diagram of complexes:

$$
\xymatrix{
CT_{0}D\ar@{^{(}->}[rr]^{\sim }\ar@{^{(}->}[dr]&&CT_{*}D\ar@{->>}[dl]\\
&CTD&
}
\eqno(\numb\label{eq4CD})
$$

The upper arrow is a quasi-isomorphism.

\bigskip

Note that precisely this proposition shows that the first term of Vassiliev spectral sequence is bounded by the upper
line $q=-(d-1)p-1$. Really, $p=-i$, $q=di-j$. For a fixed $i$ the minimal $j$ equals $i+1$ in $CT_0D$ (note, in
$CT_*D$ the minimal $j$ is $[\frac{i+1}2]$). Therefore the maximal $q$ is $di-(i+1)=(d-1)i-1=-(d-1)p-1$.

\bigskip

The complexes $CT_*D$, $CTD$, $CT_0D$ are direct sums over $i$ of finite complexes spanned by  diagrams of
complexity $i$ --- the grading $j$ is bounded: $i+1\leq j\leq 2i$ for the complexes $CTD$ and $CT_0D$;
$[\frac{i+1}2]\leq j\leq 2i$ for $CT_*D$. Let $\chi_i$, $\chi^*_i$, $\chi_i^0$ denote Euler characteristics of the
complexity $i$ components of these complexes. Consider the generating functions:
$$
\chi(t)=\sum_{i=0}^{+\infty}\chi_it^i,\qquad
\chi_*(t)=\sum_{i=0}^{+\infty}\chi_i^*t^i,\qquad
\chi_0(t)=\sum_{i=0}^{+\infty}\chi_i^0t^i.
$$

\begin{proposition}\label{p52} {\rm \cite{T2}}
$$
\chi(t)=\frac1{1-t};\qquad
\chi_*(t)=\chi_0(t)=\frac1{1-t^2}.  \nboxm
$$
\end{proposition}

For  complexity $i=0$ we count the diagram without any point --- the trivial diagram. This diagram is the lost unity (via
the Alexander duality) of the cohomology algebra.

In Appendix~\ref{A} we describe the bases of these complexes for small complexities $i$.

\section{Differential Hopf algebra structure on the complexes $CT_*D$, $CTD$, $CT_0D$}\label{s6}

Let us define multiplication and comultiplication  on complex $CT_*D$ that together with the differential will
define differential Hopf algebra structure on it. The complexes $CTD$, $CT_0D$ will inherit this structure from
$CT_*D$. These differential Hopf algebras will be called differential Hopf algebras of $T_*$-diagrams,
$T$-diagrams and $T_0$-diagrams or simply $DHAT_*D$, $DHATD$, $DHAT_0D$. The morphisms~\ref{eq4CD} will respect this
structure:

$$
\xymatrix{
DHAT_{0}D\ar@{^{(}->}[rr]^{\sim }\ar@{^{(}->}[dr]&&DHAT_{*}D\ar@{->>}[dl]\\
&DHATD&
}
\eqno(\numb\label{eq60})
$$

\begin{remark}\label{r61}
{\rm Multiplication, comultiplication and differential define only differential bialgebra structure, but our
bialgebras are connected, therefore the antipode always exists. The formula for the antipode is given in~\cite{T2}.} \nbox
\end{remark}

In the sequel we will not make difference for the terms "bialgebra" and "Hopf algebra".

\subsection{Multiplication}

The product of two $T_*$-diagrams $D_1$ and $D_2$ is defined as the shuffle of their points on the line. Orienting monomial for each
summand is the product of the orienting monomial of $D_1$ and that of $D_2$.

\begin{example}\label{ex62}
{\rm Let $d$ be odd. Consider the product of the diagrams 
\unitlength=0.2em
\begin{picture}(15,5)
 \put(0,0){\line(1,0){15}}
 \qbezier(4,0)(7.5,6)(11,0)
\end{picture}
 and 
\begin{picture}(8,3)
 \put(0,1){\line(1,0){8}}
 \put(3,-0.2){$*$}
\end{picture}%
. An orienting monomial for the first diagram
is $\beta_1\beta_{2}\alpha_{12}$, for the second one --- $\beta_1\gamma_1^*$. Their shuffle product is:}

\begin{center}
 \unitlength=0.2em
 \begin{picture}(120,15)
  \put(0,0){
   \begin{picture}(30,15)
     \put(0,0){\line(1,0){30}}
      \qbezier(5,0)(10,9)(15,0)
      \put(24,-1.2){$*$}
      \put(4,10){\footnotesize $\beta_{1}\beta_{2}\alpha_{12}\beta_{3}\gamma_{3}^*$}
   \end{picture}
   } 
  
  \put(36.5,0){$+$}
  
  \put(45,0){
   \begin{picture}(30,15)
     \put(0,0){\line(1,0){30}}
      \qbezier(5,0)(15,13)(25,0)
      \put(14,-1.2){$*$}
      \put(4,10){\footnotesize $\beta_{1}\beta_{3}\alpha_{13}\beta_{2}\gamma_{2}^*$}
   \end{picture}
   } 
   
   \put(81.5,0){$+$}
   
  \put(90,0){
   \begin{picture}(30,15)
     \put(0,0){\line(1,0){30}}
      \qbezier(15,0)(20,9)(25,0)
      \put(4,-1.2){$*$}
      \put(4,10){\footnotesize $\beta_{2}\beta_{3}\alpha_{23}\beta_{1}\gamma_{1}^*$}      
   \end{picture}
   } 
   
 \end{picture}\,\, .\nbox
\end{center}

\end{example}

Note that the trivial diagram 
\lq\lq\TD\rq\rq
 \, is the unity for this algebra.

\subsection{Comultiplication}
Comultiplication for these complexes is the map dual to concatenation (one can call it  {\it coconcatenation}).
Consider a $T_*$-diagram $D$. We say that a point $t\in\R^1$ is {\it separating} for $D$ if

1) it is not a vertex of $D$;

2) $D$ does not have edges whose left vertex belongs to $(-\infty,t)$ and the right one --- to $(t,+\infty)$.

\medskip

Two separating points are said to be {\it equivalent} if there is no vertices of $D$ between them.

\smallskip

Any separating point split $D$ into two diagrams: the left subdiagram $L_t(D)$ that is on the left from $t$; and
the right subdiagram $R_t(D)$ that is on the right from $t$. The orienting monomials of $L_t(D)$ (resp. $R_t(D)$)
are obtained from the orienting monomial of $D$ by removing the orienting elements corresponding to the right
(resp. left) subdiagram.

We define the coproduct on $D$ as follows:

$$
\Delta D=\sum_{s\in{\cal S}(D)}(-1)^{\varepsilon_s}L_s(D)\otimes R_s(D), \eqno(\numb)\label{eq63}
$$
where ${\cal S}(D)$ is the set of equivalence classes of separating points for $D$. By abuse of the language
$L_s(D)$, $R_s(D)$ denote $L_t(D)$, $R_t(D)$ for any point $t\in s$. The sign $(-1)^{\varepsilon_s}$ is obtained
as follows: Consider the product of the orienting monomial of $L_s(D)$ to the orienting monomial of $R_s(D)$. This
product differs from the orienting monomial of $D$ by some permutation  of the orienting elements.
$(-1)^{\varepsilon_s}$ is the sign of this permutation.

\medskip

For example let $1$ denote the trivial diagram \lq\lq\TD\rq\rq
 \, , then $\Delta 1=1\otimes 1$. For all the other
diagrams coproduct $\Delta D$ has at least two summands $D\otimes 1$ and $1\otimes D$.

\begin{remark}\label{r64}
{\rm The multiplication is graded commutative, but the comultiplication is not graded cocommutative. In~\cite{T2} it was
proven that the commultiplication is graded cocommutative on the homology level for any field of coefficients for
$DHATD^{odd}$, $DHATD^{even}$ and $DHAT_*D^{even}$. In the case of $DHAT_*D^{odd}$, it was proven that this is so  
over $\Q$. Theorem~\ref{t116} implies that the homology bialgebra of $DHAT_*D^{odd}$ is  graded bicommutative for any 
field of coefficients.}  \nbox
\end{remark}

\part{Complexes dual to $CTD$, $CT_*D^{even}$, $CT_0D$. Hochschild complexes}\label{s8}
In this part we describe the complexes dual to $CTD^{odd}$, $CTD^{even}$, $CT_*D^{even}$ that are the {\it
normalized Hochschild complexes of the Poisson algebras operad}, resp. {\it of the Gerstenhaber algebras operad},
resp.{\it of the Batalin-Vilkovisky algebras operad}. The dual of complex $CT_*D^{odd}$ can not be described
in terms of the Hochschild complex of none operad. Its description is given in~\cite{T1,T2,T4}.  
We describe also the duals of $CT_0D^{odd}$, $CT_0D^{even}$.

We do not prove that the constructed complexes are dual to $CTD$, $CT_*D^{even}$, $CT_0D$. However in the first
section we give a hint why this is so --- namely we describe the spaces $B_M^+$, $B_M^-$ dual to $T_M^+$,
$T_M^-$; and we show how this duality is organized.

\section{$\Z$-modules $B_M^+$, $B_M^-$}\label{ss81}
First thing to do if we want to describe the dual complexes is to describe the $\Z$-modules dual to $T_M^+$,
$T_M^-$.

Let us define the dual of $T_M^+$. Consider the $\Z$-module spanned by the monomials with \underline{commuting} generators
\underline{$\alpha_{ab}=-\alpha_{ba}$}, $a,b\in M$, $a\neq b$, that are assigned to trees with vertices labelled one-to-one by
$M$ (but not quotiented by 3-term relations~\eqref{3T}), see Section~\ref{s21}. 
This $\Z$-module is self-dual: pairing of two monomials assigned to different trees is
zero; pairing of two monomials corresponding to the same tree is equal to 1, if they define the same orientation
of the tree, and to $-1$ if they define opposite orientation. The module $T_M^+$ is a quotient-space of the
described space. Hence the dual of $T_M^+$ is a subspace of that space determined by {\it 3-term equations}. This
dual is the orthogonal to the space spanned by left-hand sides of Figure~\ref{fig2}.

The dual of $T_M^-$ is defined  in the same way, except that we take \underline{anticommuting} generators
\underline{$\alpha_{ab}=\alpha_{ba}$}, $a,b\in M$, $a\neq b$.

This definition of the duals of $T_M^+$, $T_M^-$ is not very easy to manipulate. Fortunately, these duals have
another and much simpler description.

Let $Lie^+(M)$ (resp. $Lie^-(M)$) be a usual free Lie algebra (resp. free Lie super-algebra with odd bracket) with
generators (resp. with even generators) $x_a$, $a\in M$. We define $B_M^+$ (resp. $B_M^-$) as its subspace linearly
spanned by the brackets containing each generator exactly once.

\begin{example}\label{ex81}
{\rm $B_{\{1\}}^+$, $B_{\{1\}}^-$ are spanned by the only element $x_1$.

$B_{\{1,2\}}^+$, $B_{\{1,2\}}^-$ are spanned by $[x_1,x_2]$, $[x_2,x_1]$ that are equal up to a sign.

$B_{\{1,2,3\}}^+$, $B_{\{1,2,3\}}^-$ are spanned by $[[x_1,x_2],x_3]$, $[[x_1,x_3],x_2]$, $[[x_3,x_1],x_2]$,
$[x_2,[x_1,x_3]]$, \etc. \nbox }
\end{example}

In the case $M=\{1,2,\dots,n\}$ we will write $B_n^+$, $B_n^-$ instead of $B_M^+$, $B_M^-$.

\begin{remark}\label{r82}
{\rm $B_n^+$ is the $n$-th component $\Lie(n)$ of the Lie algebras operad $\Lie$. It is well known that $\Lie(n)$ 
is isomorphic to $\Z^{(n-1)!}$. $B_n^-$ is the $n$-th component
of the operad of Lie super-algebras with odd bracket. \nbox}
\end{remark}

\begin{lemma}\label{l83}
There is a natural isomorphism of $S_n$-modules
$
B_n^+\simeq B_n^-\otimes (\pm\Z). \nboxm
$
\end{lemma}

\begin{statement}\label{st85}
The space $B_M^+$ (resp. $B_M^-$) is the dual $\Z$-module of $T_M^+$ (resp. $T_M^-$). In the case
$M=\{1,2,\dots,n\}$ this duality is $S_n$-equivariant. \nbox
\end{statement}

\noindent {\bf Sketch of the proof:} Let us define a map $\Psi$ from $B_M^+$ (resp. $B_M^-$) to the described dual
of $T_M^+$ (resp. $T_M^-$). We will define $\Psi$ inductively. If $M=\{a\}$ (cardinality of $M$ equals 1) we set
$\Psi(x_a)$ to be the only trivial tree with one vertex $a$.  In terms of monomials this is simply 1. Let cardinality $\#M$ be
$\geq 2$. Consider any bracket $L\in B_M^+$ (resp. $L\in B_M^-$). Obviously, $L=[A,B]$, where $A$, $B$ are some
subbrackets. The set $M$ is split into two subsets $M_A$ and $M_B$, where $M_A$, $M_B$ are elements of $M$
corresponding to generators in $A$, resp. in $B$. We define  $\Psi(L)=\Psi([A,B])$ by the following formula:
$$
\Psi([A,B])=\Psi(A)\Bigl(\sum\limits_{a\in M_A\atop b\in M_B}\alpha_{ab}\Bigr)\Psi(B).
\eqno(\numb)\label{eq86}
$$
It can be easily verified that

1) $\Psi$ respects (super)anticommutativity and (super)-Jacobi identities;

2) The image of $\Psi$ belongs to the dual of $T_M^+$ (resp. $T_M^-$).

\medskip

To see that $\Psi$ is bijective one needs to find a basis in $B_M^+$, $B_M^-$ that is dual to the basis of
alternated trees.

\begin{definition}\label{d87}
{\rm Let $M$ be an ordered set. A bracket in $B_M^+$ or $B_M^-$ is said to be {\it monotone} if for any its
subbracket (including itself) the generator with the minimal index occupies the left-most position and the
generator with the maximal index occupies the right-most position. \nbox }
\end{definition}

For example $[[x_1,x_3],[x_2,x_4]]$ is a monotone bracket.

\medskip

$\Psi$ of any bracket is a sum of trees.
For any monotone bracket in this sum there is  only one alternated tree. For instance, for the bracket
$[[x_1,x_3],[x_2,x_4]]$ this tree is \THREE. Hence the images of monotone brackets form a basis in the spaces dual
to $T_M^+$, $T_M^-$, and therefore the monotone brackets form a basis in $B_M^+$, $B_M^-$. Of course, this basis is
dual to the basis of alternated trees in $T_M^+$, $T_M^-$. \nbox

\section{Lie, pre-Lie and brace algebra structures on the space of any graded linear operad}\label{ss82}
In this and the next sections we briefly remind algebraic structures that exist on a graded linear operad, \cf~\cite{GV}. In
Section~\ref{ss84} we apply these structures to describe the complexes dual to $CTD$, $CT_*D^{even}$. For
people not-familiarized with the notion of operad we recommend to read~\cite{T4} instead of reading these three sections.

Let $\O=\{\O(n),n\ge 0\}$ be a graded linear operad. By abuse of the language the space $\bigoplus_{n\ge 0}\O(n)$
will be also denoted by $\O$. A tilde over an element will always designate its grading. For any element
$x\in\O(n)$ we put $n_x:=n-1$. The numbers $n$ and 1 here correspond to $n$ inputs and to 1 output respectively.

Define a new grading $|\,.\,|$ on the space $\O$. For an element $x\in\O(n)$ we put $|x|:=\tilde x+n_x=\tilde x
+n-1$. It turns out that $\O$ is a graded Lie algebra with respect to the grading $|\,.\,|$. Note that the
composition operations respect this grading.

Define the following collection of multilinear operations on the space $\O$.
$$
x\{x_1,\dots,x_n\}:=\sum(-1)^\epsilon x(id,\dots,id,x_1,id,\dots,id, x_n,id,\dots,id) \eqno(\numb)\label{eq88}
$$
for $x,x_1,\dots,x_n\in\O$, where the summation runs over all possible substitutions of $x_1,\dots,x_n$ into $x$
in the prescribed order, $\epsilon:={\sum_{p=1}^nn_{x_p}r_p}+{n_x\sum_{p=1}^n\tilde x_p}+ {\sum_{p<q}n_{x_p}\tilde
x_q}$, $r_p$ being the total number of inputs in $x$ going after $x_p$. For instance, for $x\in \O(2)$ and
arbitrary $x_1,x_2\in\O$
$$
x\{x_1,x_2\}=(-1)^{n_{x_1}+(\tilde x_1+\tilde x_2)+n_{x_1}\tilde x_2} x(x_1,x_2).
$$

By convention:
$$
x\{\}:=x.\eqno(\numb)\label{eq88'}
$$

One can check immediately the following identities:
$$
x\{x_1,\dots,x_m\}\{y_1,\dots,y_n\}=
$$
$$
\sum_{0\le i_1\le j_1\le \dots \le i_m\le j_m\le n}(-1)^{\epsilon}
x\{y_1,\dots,y_{i_1},x_1\{y_{i_1+1},\dots,y_{j_1} \},y_{j_1+1},\dots,y_{i_m},
$$
\nopagebreak
$$
x_m\{y_{i_m+1},\dots,y_{j_m}\},y_{j_m+1},\dots,y_n\}, \eqno(\numb)\label{eq89}
$$
where $\epsilon=\sum_{p=1}^m\bigl( |x_p|\sum_{q=1}^{i_p}|y_q|\bigr)$.

\begin{definition}\label{d811}
{\rm \cite{GV} A {\it brace algebra} is a graded linear space endowed with a set of multilinear $(n+1)$-ary operations
$b_n(x,x_1,\dots,x_n)=x\{x_1,\dots,x_n\}$, $n=0,1,2,\dots$, respecting the grading $|\, . \, |$ and satisfying~\eqref{eq88'},~\eqref{eq89}. \nbox }
\end{definition}

Define a bilinear operation  (respecting the grading $|\, .\,|$) $\circ$ on the space $\O$:
$$
x\circ y:=x\{y\}, \eqno(\numb)\label{eq812}
$$
for $x,y\in\O$. This operation is not associative.

\begin{definition}\label{d813}
{\rm A graded vector space $A$ with a bilinear operation
$$
\circ:A\otimes A\to A
$$
is called a {\it pre-Lie algebra}, if for any $x,y,z\in A$ the following holds:
$$
(x\circ y)\circ z -x\circ (y\circ z)=(-1)^{|y||z|}((x\circ z)\circ y- x\circ (z\circ y)).\,\,\Box
$$
}
\end{definition}

Any graded pre-Lie algebra $A$ can be considered as a graded Lie algebra with the bracket
$$
[x,y]:=x\circ y-(-1)^{|x||y|}y\circ x. \eqno(\numb\label{bracket})
$$

The description of the operad of pre-Lie algebras is given in~\cite{Cha}.

The following lemma is a corollary of the identity~(\ref{eq89}) applied to the case $m=n=1$.

\begin{lemma}\label{815}
The operation~\eqref{eq812} defines a graded pre-Lie algebra structure on the space $\O$.~$\Box$
\end{lemma}

In particular this lemma implies that any graded linear operad $\O$ can be considered as a graded Lie algebra with
the bracket~\eqref{bracket}. This bracket is usually called {\it Gerstenhaber bracket} in honor of Murray Gerstenhaber who 
discovered this operation for the Hochschild cochain complex of an associative algebra, see~\cite{Ge} and also next section.

\section{Hochschild complexes}\label{ss83}

Let $\O =\bigoplus_{n\ge0}\O(n)$ be a graded linear operad equipped with a morphism
$$
\Pi:\Assoc\to\O
$$
from the operad $\Assoc$. This morphism defines the element $m=\Pi (m_2)\in\O(2)$, where the element
$m_2=x_1x_2\in\Assoc(2)$ is the operation of multiplication. Note that the elements $m_2$, $m$ are odd with
respect to the new grading $|\,.\,|$ ($|m|=|m_2|=1$) and $[m,m]=[\Pi(m_2),\Pi(m_2)]=2\Pi(m_2\circ m_2)=0$. (One
has $m_2\circ m_2=-(x_1\cdot x_2)\cdot x_3 +x_1\cdot (x_2\cdot x_3)= -x_1\cdot x_2\cdot x_3+x_1\cdot x_2\cdot
x_3=0$.) Thus $\O$ becomes a differential graded Lie algebra with the differential $\partial$:
$$
\partial x:=[m,x]=m\circ x-(-1)^{|x|}x\circ m,
\eqno(\numb)\label{eq816}
$$
for $x\in \O$.

We will call complex $(\O,\partial)$    {\it Hochschild complex} of  operad $\O$. Actually a better name
would be {\it Hochschild complex} of the {\it morphism} $\Pi:\Assoc\to\O$ since this complex is in fact the
deformation complex of the morphism $\Pi$, see~\cite{K2,KS}. We preferred the first name for its shortness.

\begin{example}\label{ex817}
{\rm If $\O$ is the endomorphism operad $\Endom(A)$ of a vector space $A$, and we have a morphism
$$\Pi:\Assoc\to \Endom(A),
$$
that defines an associative algebra structure on $A$, then the corresponding complex
$\bigl(\bigoplus_{n=0}^{+\infty} Hom(A^{\otimes n},A),
\partial\bigr)$ is the usual Hochschild cochain complex
$C^*(A,A)$ of an associative algebra $A$.~$\Box$ }
\end{example}
%

Define another grading
$$
deg:=|\,.\,|+1 \eqno(\numb)\label{eq818}
$$
on the space $\O$. With respect to this grading the bracket $[\, .\, ,\, .\, ]$ is homogeneous of degree $-1$.

It is easy to see that the product $*$, defined as follows
$$
x*y:=(-1)^{|x|}m\{x,y\}=(-1)^{\tilde y(n_x+1)}m(x,y), \eqno(\numb)\label{eq819}
$$
for $x,y\in\O$, together with the differential $\partial$ defines a differential graded associative algebra
structure on $\O$ with respect to the grading~\eqref{eq818}.

\begin{theorem}\label{819}
{\rm \cite{Ge,GV}} The multiplication $*$ and the bracket $[\, .\, ,\, .\, ]$ induce a Gerstenhaber algebra structure on the
homology of the Hochschild complex $(\O,
\partial)$.~$\Box$
\end{theorem}

For a definition of Gerstenhaber algebras we refer to the next section.

\medskip

\noindent{\bf Proof:} The proof is deduced from the following homotopy formulas.
$$
x*y-(-1)^{deg(x)deg(y)}y*x=(-1)^{deg(x)}(\partial(x\circ y)-
\partial x\circ y-(-1)^{deg(x)-1}x\circ \partial y).
\eqno(\numb)\label{eq820}
$$
The above formula proves the graded commutativity of the multiplication~$*$.

\vspace{2mm}

\noindent $[x,y*z]-[x,y]*z-(-1)^{(deg(x)-1)deg(y)}y*[x,z]=$
$$
\mbox{\footnotesize$=(-1)^{deg(x)+deg(y)} \bigl(\partial(x\{y,z\})-(\partial x)\{y,z\}-(-1)^{|x|}x\{\partial y,
z\}-(-1)^{|x|+|y|}x\{y,\partial z\}).$} \eqno(\numb)\label{eq821}
$$
This formula proves the compatibility of the bracket with the multiplication.~$\Box$

\section{Dual complexes}\label{ss84}

\begin{definition}\label{d822}
{\rm A graded commutative algebra is called {\it Poisson algebra} if it is endowed with a Lie bracket respecting
the grading and  compatible with the multiplication:
$$
[x,yz]=[x,y]z+(-1)^{\tilde{x}\tilde{y}}y[x,z].~\square \eqno(\numb)\label{eq823}
$$
}
\end{definition}

\begin{definition}\label{d824}
{\rm A graded commutative algebra is called {\it Gerstenhaber algebra} if it is endowed with a Lie bracket of
degree $-1$ and  compatible with the multiplication:
$$
[x,yz]=[x,y]z+(-1)^{(\tilde{x}-1)\tilde{y}}y[x,z].~\square \eqno(\numb)\label{eq825}
$$
}
\end{definition}

\begin{definition}\label{d826}
{\rm A Gerstenhaber algebra is called {\it Batalin-Vilkovisky algebra} if it is endowed with an unary linear
operation $\delta$ of degree $-1$ satisfying:

(i) $\delta^2=0$;

(ii) $\delta(ab)=\delta(a)b+(-1)^{\tilde a}a\delta(b)+(-1)^{\tilde a} [a,b]. \nboxm$ }
\end{definition}

Note that (i) and (ii) imply

(iii) $\delta([a,b])=[\delta(a),b]+(-1)^{\tilde a +1}[a,\delta(b)].$

\bigskip

Denote by $\Poiss$, $\Gerst$, $\BV$ the operads of Poisson, Gersenhaber and Batalin-Vilkovisky algebras.

Let us describe these operads.

The space $\Poiss(n)$ of all multilinear $n$-ary operations that are induced by a Poisson algebra structure is
described as follows. Consider a free Lie algebra $Lie(x_1,\dots,x_n)$ with $n$ generators $x_1,\dots,x_n$ and
consider the symmetric algebra $S^*Lie(x_1,\dots,x_n)$. Any symmetric algebra of a Lie algebra is endowed with a
Poisson algebra structure. This one is a {\it free Poisson algebra}. The space $\Poiss(n)$ is a subspace of this
algebra linearly spanned by the products of brackets using each generator exactly once. For instance, if $n=3$ we
have the following elements: $x_1\cdot x_2\cdot x_3$, $x_1\cdot [x_2,x_3]$, $[[x_3,x_1],x_2]$, \etc. The space $\Poiss$ is bigraded:
the first bigrading {\it complexity} $i$ is the total number of comas in the products of brackets; the second
grading $j=n$ is the component number.

The space $\Gerst(n)$ is defined analogously. The only difference is that we need to consider a free graded Lie
algebra $Lie_1(x_1,\dots,x_n)$ with the bracket of degree $-1$ but also with the generators $x_1,\dots,x_n$ of
degree zero. The space $\Gerst$ is also bigraded by complexity and by  component number.

To define $\BV(n)$ we need to start from a free graded Lie algebra $Lie_1(x_1,\dots,x_n,\delta x_1,\dots,\delta
x_n)$ with the bracket of degree $-1$ and with the generators $x_1,\dots,x_n$ of degree $0$ and $\delta
x_1,\dots,\delta x_n$ of degree $-1$. Then we consider the symmetric algebra $S^*Lie_1(x_1,\dots,x_n,\delta
x_1,\dots,\delta x_n)$ and take its subspace linearly spanned by the products of brackets using each index
$1,\dots,n$ exactly once. For instance, if $n=3$ we have the elements $x_1\cdot \delta x_2\cdot x_3$,
$x_1\cdot [x_2,\delta x_3]$, $[[\delta x_3,x_1],\delta x_2]$, \etc. The grading complexity is the total number of comas
and deltas. The second grading is as usual the component number.

\bigskip

Any Poisson (resp. Gerstenhaber, resp. Batalin-Vilkovisky) algebra is a commutative algebra and therefore is an
associative algebra. It means that the operads $\Poiss$, $\Gerst$, $\BV$ are endowed with a map from the operad
$\Assoc$ of associative algebras. By the previous subsection they form Hochschild complexes that we will denote by
$(\Poiss,\partial)$, $(\Gerst,\partial)$, $(\BV,\partial)$. These are not yet the duals of $CTD^{odd}$,
$CTD^{even}$, $CT_*D^{even}$.

There are two types of products of brackets that span $\Poiss$, $\Gerst$ or $\BV$: in the first group we put
those products that have at least one  factor $x_i$ for some $i$; in the second group --- all the others. The point is that both
groups span a subcomplex (and therefore each of the complexes $(\Poiss,\partial)$, $(\Gerst,\partial)$,
$(\BV,\partial)$ is a direct sum of two complexes). It is easy to show that the first subcomplex is always acyclic
and even contractible. This implies that the second one is quasi-isomorphic to $(\Poiss,\partial)$,
$(\Gerst,\partial)$, $(\BV,\partial)$ respectively. The second subcomplex  will be called  {\it normalized
Hochschild complex} and denoted by $(\Poiss^{Norm},\partial)$, $(\Gerst^{Norm},\partial)$,
$(\BV^{Norm},\partial)$.

\begin{theorem}\label{t827}
{\rm \cite{T2}} The complexes $(\Poiss^{Norm},\partial)$, $(\Gerst^{Norm},\partial)$, $(\BV^{Norm},\partial)$ are dual
to $CTD^{odd}$, $CTD^{even}$, $CT_*D^{even}$. \nbox
\end{theorem}

\noindent {\bf Idea of the proof:} It is a direct check.  To prove one needs to use the duality described in Secion~\ref{ss81}. \nbox

\bigskip
 
These complexes are differential Hopf algebras (since they are dual to differential Hopf algebras), they also 
inherit  pre-Lie, Lie and brace algebra structures.

\bigskip

To define the complexes dual to $CT_0D^{odd}$, $CT_{0}D^{even}$ one needs to quotient the complexes
$(\Poiss^{Norm},\partial)$, $(\Gerst^{Norm},\partial)$ by the {\it neighbor commutativity relations}. Namely,
instead of taking free (graded) Lie algebras $Lie(x_1,\dots,x_n)$ or $Lie_1(x_1,\dots,x_n)$  one should take Lie
algebras with the same generators but satisfying  the relations $[x_i,x_{i+1}]=0$, for $i=1,\dots,n-1$. And then we
proceed as before. We denote these complexes by $(\Poiss^{zero},\partial)$, $(\Gerst^{zero},\partial)$. These
complexes are as before differential Hopf algebras, but they have no more pre-Lie, Lie or brace algebra
structure.

\section{Bases of monotone bracket diagrams}\label{ss85}
In Appendix~\ref{A} we describe bases of the complexes $(\Poiss^{Norm},\partial)$, $(\Gerst^{Norm},\partial)$,
$(\BV^{Norm},\partial)$,
$(\Poiss^{zero},\partial)$,
$(\Gerst^{zero},\partial)$
 for small complexities $i$. This provides an illustration of the
given below considerations.

\medskip

Products of brackets that span the complexes $(\Poiss^{Norm},\partial)$, $(\Gerst^{Norm},\partial)$ will be called
{\it bracket diagrams}. If none of the factors of a bracket diagram contains subbracket $[x_i,x_{i+1}]$ for some
$i$, then this bracket will be called {\it bracket zero-diagram}. For example,  $[x_1,x_3]\cdot [x_2,[x_4,x_5]]$
is not bracket zero-diagram because it contains subbracket $[x_4,x_5]$. If all the factors of a bracket
diagram are monotone brackets, see Section~\ref{ss81}, then this diagram is called {\it monotone
bracket diagram}. Monotone bracket diagrams (resp. monotone bracket zero-diagrams) form  bases in the complexes
$(\Poiss^{Norm},\partial)$, $(\Gerst^{Norm},\partial)$ (resp. $(\Poiss^{zero},\partial)$,
$(\Gerst^{zero},\partial)$) that are dual to the bases of alternated $T$-diagrams (resp. $T_0$-diagrams).

\medskip

Let us make a few remarks about how these bases behave with respect to the algebraic operations in the
corresponding complexes.

Note first that the bases of alternated $T/T_*/T_0$-diagrams are invariant with respect to multiplication and
comultiplication in the sense that all the summands in the shuffle product and in the formula~\eqref{eq63} of
comultiplication will be basis elements (if we apply multiplication or comultiplication to basis elements). It
means that we do not need to use basis decomposition to compute these operations. On the contrary we do need to use basis
decomposition to compute the differential of these complexes.
\medskip

We have the same for the dual complexes: multiplication and comultiplication are invariant with respect to the
bases of monotone bracket (zero)-diagrams, but the differential is not.

Neither pre-Lie product, nor brace algebra operations~\eqref{eq88} respect the basis. However in one case these
operations $x\circ y$, $x\{x_1,\dots,x_n\}$ do respect the basis. Namely, if $y$, resp. $x_1,\dots,x_n$ consist of
only one minimal component --- only one bracket. Really, in this case all the summands of~\eqref{eq88} will be
also monotone bracket diagrams. We will use this fact in Sections~\ref{s9} and~\ref{s11}.

\part{Upper diagonal of the Vassiliev spectral sequence and of the related complexes}\label{upper_diagonal}
In Section~\ref{s7} we prove that the upper diagonal homology groups of complex $CT_{0}D$ are trivial (Theorem~\ref{t71}).
Theorem~\ref{t72} describes the upper diagonal homology groups of complexes $CTD^{odd}$, $CTD^{even}$. This theorem
is not completely proven in this section.

In Section~\ref{s9} we introduce new homological operations on Hochschild complexes. These operations are defined in finite characteristic.
It will be proven elsewhere that they are related to Dyer-Lashof operations (see Theorem~\ref{cone3}). By means of these operations we find cycles in the dual complexes
and this completes the proof of Theorem~\ref{t72}.

In Section~\ref{s10} we prove some composition formulas for the operations defined in previous section.

By Theorem~\ref{t72} any upper diagonal diagram of $CTD$ is a cycle homologous to standard diagram~\eqref{eq72} with some coefficient.
In Section~\ref{s11} we use the duality of Section~\ref{ss81} and also the  formulas for dual cycles (from Section~\ref{s9}) of dual complexes
in order to determine the above coefficients. 

In Section~\ref{s14} we describe the relation between homology bialgebra of $CTD$  and that of $CT_{0}D$. Main Theorem~\ref{t116} is given 
without proof.

\section{Upper diagonal of the Vassiliev spectral sequence}\label{s7}
We have seen that for a fixed complexity $i$ the minimal $j$ equals $i+1$: in this case the diagrams have only 
one minimal component. The homology groups of the complex
$CT_0D$ in the bigradings $(i,i+1)$ are the upper diagonal elements of the first term of the Vassiliev spectral sequence.

\begin{theorem}\label{t71}
For all $i\in\N$ the homology groups in the bigradings $(i,i+1)$ are trivial for complex $CT_0D$, and are cyclic
for complex $CTD$: if the corresponding group is untrivial, as a generator one can take the one diagram cycle:
$$
\unitlength=0.3em
\begin{picture}(35,13)(0,-3)
  \put(0,0){\line(1,0){35}}
  \qbezier(5,0)(7.5,4)(10,0)
  \qbezier(5,0)(10,7)(15,0) 
  \qbezier(5,0)(15,12)(25,0)
  \qbezier(5,0)(17.5,15)(30,0)
  \put(16,1.2){\Large $\ldots$}
  \put(8,0){ $\underbrace{\hspace{6.5em}}_{\text{$i$ {\rm chords}}}$ }
\end{picture}\quad . \nboxm
\eqno(\numb)\label{eq72}
$$
\end{theorem}

\noindent {\bf Proof:} We will consider both complexes $CT_0D$ and $CTD$ simultaneously.

Since there is no diagrams of complexity $i$ with less then $i+1$ points, any diagram of bigrading $(i,i+1)$
always defines a cycle. We will prove that any such cycle is homologous to diagram~\eqref{eq72} with some
coefficient\footnote{In Section~\ref{s11} we determine these coefficients.}. But diagram~\eqref{eq72}  
does not belong to the space of $T_0$-diagrams, therefore this coefficient
must be zero if we treat complex $CT_0D$.

We will say that an alternated $T$-diagram (resp. $T_0$-diagram) is of {\it $r$-type}, where $2\leq r\leq i+1$, if
it has the edges $(1,r)$, $(1,r+1)$,\dots, $(1,i+1)$ and for the points $r+1$, $r+2$,\dots, $i+1$ there is no any
other incident edges.

\begin{figure}[!h]
\begin{center}
\unitlength=0.15em
\begin{picture}(50,23)(-3,0)
 \put(-3,0){\line(1,0){56}}
 \qbezier(0,0)(25,30)(50,0)
 \qbezier(0,0)(20,20)(40,0)
 \qbezier(0,0)(15,15)(30,0)
 \qbezier(0,0)(10,10)(20,0)
 \qbezier(10,0)(20,10)(30,0)
 \put(-1.5,-5){$1$}
 \put(8.5,-5){$2$}
 \put(18.5,-5){$3$}
 \put(28.5,-5){$4$}
 \put(38.5,-5){$5$}
 \put(48.5,-5){$6$}
\end{picture}
\end{center}

\caption{diagram $D$ of 4-type. $i=5$.}\label{fig5''}
\end{figure}
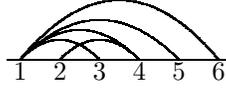

Any alternated diagram of bigrading $(i,i+1)$ is of $(i+1)$-type. The diagram~\eqref{eq72} is the only one of
2-type.

Let us prove that any alternated diagram $D$ of $r$-type ($r\ge 3$) is homologous to a sum of alternated diagrams
of $(r-1)$-type. Consider a diagram $D_\bigstar$ of bigrading $(i,i+2)$ that is obtained from $D$ by splitting the
point $r$ into two points: left one $r_-$ and right one $r_+$. The edges that were not incident to $r$ stay
without any change. The edge $(1,r)$ becomes $(1,r_+)$; all the other edges incident to $r$ become incident to
$r_-$ (if there were no such edges the diagram $D$ was already of type $r-1$), see the example below --- diagram
$D_\bigstar$ is obtained from the one of the Figure~\ref{fig5''}.

\begin{figure}[!ht]
\begin{center}
\unitlength=0.15em
\begin{picture}(60,27)(-3,0)
 \put(-3,0){\line(1,0){66}}
 \qbezier(0,0)(30,39)(60,0)
 \qbezier(0,0)(25,27)(50,0)
 \qbezier(0,0)(20,20)(40,0)
 \qbezier(0,0)(10,10)(20,0)
 \qbezier(10,0)(20,10)(30,0)
 \put(-1.5,-5){$1$}
 \put(8.5,-5){$2$}
 \put(18.5,-5){$3$}
 \put(27,-5){$4_-$}
 \put(37.5,-5){$4_+$}
 \put(48.5,-5){$5$}
 \put(58.5,-5){$6$}

\end{picture}
\end{center}

\caption{diagram $D_\bigstar$}\label{fig6}
\end{figure}
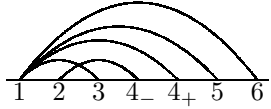

\begin{figure}[!ht]
\begin{center}
\unitlength=0.1em
\begin{picture}(270,30)
 
 \multiput(0,0)(70,0){4}{$\pm$}
 
 \put(10,0){
\begin{picture}(50,23)
 \put(-3,0){\line(1,0){56}}
 \qbezier(0,0)(25,30)(50,0)
 \qbezier(0,0)(20,20)(40,0)
 \qbezier(0,0)(15,15)(30,0)
 \qbezier(0,0)(10,10)(20,0)
 \qbezier(0,0)(5,6)(10,0)
\end{picture}
}

 \put(80,0){
\begin{picture}(50,23)
 \put(-3,0){\line(1,0){56}}
 \qbezier(0,0)(25,30)(50,0)
 \qbezier(0,0)(20,20)(40,0)
 \qbezier(0,0)(15,15)(30,0)
 \qbezier(0,0)(5,6)(10,0)
 \qbezier(10,0)(15,6)(20,0)
\end{picture}
}

 \put(150,0){
\begin{picture}(50,23)
 \put(-3,0){\line(1,0){56}}
 \qbezier(0,0)(25,30)(50,0)
 \qbezier(0,0)(20,20)(40,0)
 \qbezier(0,0)(15,15)(30,0)
 \qbezier(0,0)(10,10)(20,0)
 \qbezier(10,0)(15,6)(20,0)
\end{picture}
}

\put(220,0){
\begin{picture}(50,23)
 \put(-3,0){\line(1,0){56}}
 \qbezier(0,0)(25,30)(50,0)
 \qbezier(0,0)(20,20)(40,0)
 \qbezier(0,0)(15,15)(30,0)
 \qbezier(0,0)(10,10)(20,0)
 \qbezier(10,0)(20,10)(30,0)
\end{picture}
}

\end{picture}
\end{center}

\caption{$\partial D_\bigstar$}\label{fig7}
\end{figure}
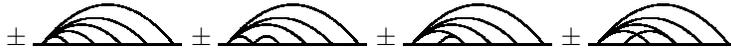

The diagram $D_\bigstar$ is also alternated $T$-diagram (resp. $T_0$-diagram).  Consider the differential of this
diagram $\partial D_\bigstar$, see Figure~\ref{fig7}.

The last admissible gluing --- gluing of the points $r_-$ and $r_+$
--- gives the diagram $D$. All the other boundary diagrams have the edges $(1,r_+)$, $(1,r+1)$, $(1,r+2)$,\dots,
$(1,i+1)$ and they are the only incident for the points $r_+$, $r+1$, $r+2$,\dots, $i+1$ respectively. For any of
these boundary diagrams its subdiagram on the first $r$ points is a tree. It means that the decomposition of the
sum of these diagrams in our basis is a sum of alternated diagrams of type $r-1$. This completes the proof of the
theorem. \nbox

\bigskip
\bigskip

\begin{theorem}\label{t72}
For complex $CTD^{even}$ the homology groups of the bigrading $(i,i+1)$ are
$$
\begin{cases}
\Z,&i=1;\\
\Z_p,&\text{$i=p^k$, where $p$ is any prime, $k\in\N$;}\\
0,&\text{otherwise.}
\end{cases}
\eqno(\numb)\label{eq73}
$$

For complex $CTD^{odd}$ these groups are
$$
\begin{cases}
\Z,&i=1,\,\, 2;\\
\Z_p,&\text{$i=2p^k$, where $p$ is any prime, $k\in\N$;}\\
0,&\text{otherwise.}
\end{cases}
\eqno(\numb)\label{eq74}
$$
$\Box$
\end{theorem}

\noindent {\bf First part of the proof:} In this section  we will prove that the homology groups are some
quotient-groups of the described groups. To finish the proof we will present the dual cycles in the dual complexes
tensored to the corresponding cyclic group~\eqref{eq73},~\eqref{eq74}, see Section~\ref{s9}.

\begin{definition}\label{d75}
{\rm A {\it divided product} $\langle D_1,D_2\rangle $ of two diagrams $D_1$ and $D_2$ is the sum of those elements in the shuffle
product $D_1*D_2$ that have the left-most point of $D_1$ on the left from the left-most point of $D_2$.} \nbox
\end{definition}

Note $D_1*D_2=\langle D_1,D_2\rangle +(-1)^{\deg D_1 \deg D_2}\langle D_2,D_1\rangle $, where $\deg D_i$ is the number of orienting elements
of $D_i$, or, what is the same modulo 2,  the corresponding cohomology degree of the knot space,   \ie
$p+q=i(d-1)+j$.

\begin{figure}[!ht]

\begin{center}
 \unitlength=0.2em
 \begin{picture}(140,15)
   \put(0,0){$\langle$}
   \put(2,0){
   \begin{picture}(20,10)
     \put(0,0){\line(1,0){20}}
      \qbezier(5,0)(10,9)(15,0)
    \end{picture}
   } 
  
\put(26.5,0){$,$}

   \put(30,0){
   \begin{picture}(15,10)
     \put(0,0){\line(1,0){15}}
      \put(6.5,-1.2){$*$}
    \end{picture}
   } 
   
   \put(48,0){$\rangle$}

\put(56.5,0){$=$}
  
  \put(65,0){
   \begin{picture}(30,10)
     \put(0,0){\line(1,0){30}}
      \qbezier(5,0)(10,9)(15,0)
      \put(24,-1.2){$*$}
    \end{picture}
   } 
  
  \put(101.5,0){$+$}
  
  \put(110,0){
   \begin{picture}(30,10)
     \put(0,0){\line(1,0){30}}
      \qbezier(5,0)(15,13)(25,0)
      \put(14,-1.2){$*$}
   \end{picture}
   }

 \end{picture}
\end{center}

\caption{Example of a divided product}\label{fig8}
\end{figure}

Denote the diagram~\eqref{eq72} by $Z_i$.

In the case of odd $d$ this diagram has $i+1$ orienting elements, each of them is a  representative of a point. We fix its orientation
by the monomial:
$$
\beta_{1}\beta_{2}\dots\beta_{i+1}.
$$ 
We orient the edges all from left to right:
$$
\unitlength=0.3em
\begin{picture}(35,13)(0,0)
  \put(0,0){\line(1,0){35}}
  \qbezier(5,0)(7.5,4)(10,0)
  \put(10,0){\vector(2,-3){0}}
  \qbezier(5,0)(10,7)(15,0) 
  \put(15,0){\vector(2,-3){0}}
  \qbezier(5,0)(15,12)(25,0)
  \put(25,0){\vector(2,-3){0}}
  \qbezier(5,0)(17.5,15)(30,0)
  \put(30,0){\vector(2,-3){0}}
  \put(16,1.2){\Large $\ldots$}
\end{picture}
\eqno(\numb)\label{eq76}
$$

\medskip

In the case of even $d$ we fix an orientation of $Z_i$ by the orienting monomial:
$$
\beta_1\alpha_{12}\beta_{2}\alpha_{13}\beta_{3}\dots\alpha_{1(i+1)}\beta_{i+1}.
$$

Consider  divided product $\langle Z_k,Z_l\rangle $, $k,l\geq 1$.

For even $d$ one has:
$$
\partial\langle Z_k,Z_l\rangle ={{k+l}\choose k} Z_{k+l}.
\eqno(\numb)\label{eq77}
$$

For odd $d$ one has
$$
\partial\langle Z_k,Z_l\rangle =(-1)^k{{k+l}\choose k}_{-1} Z_{k+l},
\eqno(\numb)\label{eq78}
$$
where 
$$
{{k+l}\choose k}_{-1}=
\begin{cases}
0,&\text{$k$ and $l$ are odd;}\\
\left({{\left[\frac{k+l}2\right]}\atop{\left[ \frac k2\right]}}\right), & \text{otherwise},
\end{cases}
$$ 
see Appendix~\ref{C}. For instance, $\partial \langle Z_1,Z_{2r}\rangle =-Z_{2r+1}$.

Really, all the gluings in $\partial\langle Z_k,Z_l\rangle $ will be compensated with each other except the gluing of the
left-most point of $Z_k$ with the left-most point of $Z_l$. Such gluing always gives $Z_{k+l}$ (with some sign).
This situation occurs ${k+l}\choose k$ times since the points $2,3,\dots,k+1$ of $Z_k$ can be arbitrarily shuffled
with the points $2,3,\dots,l+1$ of $Z_l$. In the case of even $d$ all these gluings give the same orientation of
$Z_{k+l}$. In the case of odd $d$ one needs to apply the combinatorics of shuffles of anticommuting
elements, see Appendix~\ref{C}.

The following lemma completes the first part of the proof of the theorem. 

\begin{lemma}\label{l78}
The greatest common deviser of the numbers ${i\choose
1},{i\choose 2},\dots,{i\choose{i-1}}$ {\rm (}resp. ${i\choose 1}_{-1},{i\choose 2}_{-1},\dots,{i\choose{i-1}}_{-1}${\rm )} is
equal to $0$ if $i=1$ {\rm (}resp. $i=1$, $2${\rm )}, to $p$ if $i=p^k$ {\rm (}resp. $i=2p^k${\rm )}, $p$ being any prime, and to $1$ in all
the other cases. \nbox
\end{lemma}

\section{End of the proof of Theorem~\ref{t72}}\label{s9}
We will consider Hochschild complexes $(\O,\partial)$, see Section~\ref{ss83},  over some commutative ring
$\kk$. We are mostly interested in the cases when $\kk$ is $\Z$, $\Q$ or any finite field $\Z_p$.

Let $\varphi\in\O$. We will define {\it $[n]$-operation} of $\varphi$ as follows:
$$
\varphi^{[n]}:=(\dots((\underbrace{\varphi\circ\varphi)\circ\varphi)\dots)\circ\varphi}_{\text{$n$ times}}
\eqno(\numb\label{eq91})
$$
One has $\varphi^{[1]}=\varphi$, $\varphi^{[n+1]}=\varphi^{[n]}\circ\varphi$.

Remind that we defined two gradings $deg$ and $|\, .\, |$ on the complexes $(\O,\partial)$, and $deg=|\, .\, |+1$,
see formula~\eqref{eq818}. In the cases $\O$ is $\Poiss$ or $\Gerst$ the grading $deg$ is the main grading:
modulo 2 it equals to the corresponding homology degree of the knot spaces. In our formulas we will use both
gradings. $deg(\varphi^{[n]})\equiv (n\cdot deg(\varphi) +n-1)\mod 2$; $|\varphi^{[n]}|\equiv n\cdot|\varphi|\mod
2$.

\medskip

Theorem~\ref{t72} will follow from the proposition:

\begin{proposition}\label{p92}
Let $\varphi\in\O$ and $\partial\varphi=0$.

If $deg(\varphi)$ is odd then
$$
\partial(\varphi^{[n]})=-\sum_{i=1}^{n-1} {n\choose i }\varphi^{[i]}*\varphi^{[n-i]}.
\eqno(\numb\label{eq93})
$$

If $deg(\varphi)$ is even then
$$
\partial(\varphi^{[n]})=\sum_{i=1}^{n-1} (-1)^{i-1} {n\choose i}_{-1} \varphi^{[i]}*\varphi^{[n-i]}.
\eqno(\numb\label{eq94})
$$
$\Box$
\end{proposition}

\noindent{\bf Proof of Proposition~\ref{p92}:} The proof is based on two formulas. The first one is~\eqref{eq820}.
It can be rewritten as:
$$
\partial(x\circ y)=\partial x\circ y-(-1)^{deg(x)-1}x\circ \partial y+(-1)^{deg(x)}(x*y-(-1)^{deg(x)deg(y)}y*x).
\eqno(\numb)\label{eq95}
$$
The second one is given by the following lemma:

\begin{lemma}\label{l96}
$$
(x*y)\circ z=x*(y\circ z)+(-1)^{deg(y)\cdot|z|}(x\circ z)*y. \nboxm \eqno(\numb\label{eq97})
$$
\end{lemma}

\noindent{\bf Proof of Lemma~\ref{l96}:} Due to~\eqref{eq819} and~\eqref{eq89}, one has
\begin{gather*}
(x*z)\circ y=(-1)^{|x|}m\{x,y\}\{z\}=(-1)^{|x|}m\{x,y\{z\}\}+(-1)^{|x|+|y|\cdot|z|}m\{x\{z\},y\}\\
=x*(y\circ z)+(-1)^{|z|+|y|\cdot |z|}(x\circ z)*y=x*(y\circ z)+(-1)^{deg(y)\cdot|z|}(x\circ z)*y. \nboxm
\end{gather*}

\bigskip

The formulae~\eqref{eq95},~\eqref{eq97} together with the combinatorial formula
$$
{{n+1}\choose i}={n\choose i}+ {n\choose {i-1}},
$$
and its super-analog  (see Appendix~\ref{C}):
$$
{{n+1}\choose i}_{-1}={n\choose i}_{-1}+ (-1)^{n-i+1}{n\choose {i-1}}_{-1},
$$
prove Proposition~\ref{p92}. \nbox

\bigskip
\bigskip

Proposition~\ref{p92} and Lemma~\ref{l78}
imply that $[n]$-operation is a homology operation in
the following cases:

1) $\kk=\Z$ or any other commutative ring: if $deg(\varphi)$ is even and $n=2$ --- the resulting element has odd grading $deg$.
\smallskip

2) $\kk=\Z_2$: for any grading $deg$, $n=2^k$.

\smallskip

3) $\kk=\Z_p$, $p$ being any odd prime:

if $deg(\varphi)$ is odd and $n=p^k$ --- the resulting element has odd grading $deg$;

if $deg(\varphi)$ is even and $n=2p^k$ --- the resulting element has odd grading $deg$.

\begin{remark}\label{hom}
{\rm To show that some operation is a homology operation one should prove that the
operation does not depend on a representative of a cycle: if we take a homologous element the result will be 
homologous. A priori it is not evident. In the article we do not prove this for the above operations. \nbox}
\end{remark}

\vspace{4mm}

Now, let $\O$ be $\Gerst$ or $\Poiss$ and $\varphi=[x_1,x_2]$. $\varphi$ defines a nontrivial cycle for any ring
$\kk$ of coefficients since it is the only diagram of complexity $i=1$. $deg(\varphi)$ is odd in the case of the
operad $\Gerst$ and is even in the case of the operad $\Poiss$. Suppose  we are in one of the situations where we
want to prove non-triviality of the corresponding cyclic group, see the statement of Theorem~\ref{t72}, and $\kk$
is the corresponding cyclic group. Obviously, $\varphi^{[i]}$ has the complexity $i$. On the other hand, basis
decomposition of $\varphi^{[i]}$ has the monotone diagram $[[\dots[[x_1,x_2],x_3],\dots],x_{i+1}]$ with 
coefficient $\pm 1$ (since operation $\circ$ respects the basis, see Section~\ref{ss85}). Therefore this
element provides a non-trivial pairing with the diagram~\eqref{eq72} denoted by $Z_i$ and this proves a
non-triviality of the corresponding homology group. Thus Theorem~\ref{t72} is proved. \nbox

\section{Composition of operations}\label{s10}

In the previous section we defined $[n]$-operations that are homology operations in the following cases:

1) over $\Z_2$: $[2^k]$-operations for both even and odd elements. The resulting elements are always odd.

2) over $\Z_p$: $[p^k]$-operation for odd elements and $[2p^k]$-operation for even element. The resulting elements
are also odd.

\medskip

In this section we will prove two lemmas:

\begin{lemma}\label{l105'}
If $deg(\varphi)$ is odd, then
$$\varphi^{[pn]}\equiv(\varphi^{[p]})^{[n]}\mod p$$
for any prime $p$ and $n\in\N$. \nbox
\end{lemma}

\begin{lemma}\label{l1010}
If $deg(\varphi)$ is even then
$$
\varphi^{[2n]}=(\varphi^{[2]})^{[n]}, \eqno(\numb\label{eq1011})
$$
for any $n\in\N$. \nbox
\end{lemma}

Note that equality~\eqref{eq1011} is over $\Z$.

\medskip

These lemmas imply that

1) Over $\Z_2$: $[2^k]$-operation is $k$ times applied $[2]$-operation.

2) Over $\Z_p$: $[p^k]$-operation (on odd elements) is $k$ times applied $[p]$-operation; $[2p^k]$-operation (on
even elements) is $[2]$-operation and then $k$ times applied $[p]$-operation.

\smallskip

This explains why we have powers of $p$. Note, the above assertions are true already for the chains.

\bigskip

\noindent{\bf Proof of Lemma~\ref{l105'}:} Before giving a general proof let us see what is going on in the
simplest case: $deg(\varphi)$ is odd, $p=2$ and $n=2$. We do the following computations over $\Z$ and then
compare the results:
$$
\begin{array}{rl}
{\scriptstyle \varphi^{[4]}\quad =}&
{\scriptstyle ((\varphi\circ\varphi)\circ\varphi)\circ\varphi=\varphi\{\varphi\}\{\varphi\}\{\varphi\}=
\left(\varphi\{\varphi\{\varphi\}\}+2\varphi\{\varphi,\varphi\}\right)\{\varphi\} }\\
{\scriptstyle = }&{\scriptstyle \varphi\{\varphi\{\varphi\}\}\{\varphi\}+2\varphi\{\varphi,\varphi\}\{\varphi\} }\\
{\scriptstyle = }&{\scriptstyle \varphi\{\varphi\{\varphi\{\varphi\}\}\}+\varphi\{\varphi,\varphi\{\varphi\}\}+
\varphi\{\varphi\{\varphi\},\varphi\}+2\varphi\{\varphi\{\varphi\},\varphi\}+2\varphi\{\varphi,\varphi\{\varphi\}\}+
6\varphi\{\varphi,\varphi,\varphi\} }\\
{\scriptstyle = }&{\scriptstyle \varphi\{\varphi\{\varphi\{\varphi\}\}\}+3\varphi\{\varphi,\varphi\{\varphi\}\}+
3\varphi\{\varphi\{\varphi\},\varphi\}+6\varphi\{\varphi,\varphi,\varphi\}; }
\end{array}
\eqno(\numb\label{eq102})
$$

$$
\begin{array}{rl}
{\scriptstyle (\varphi^{[2]})^{[2]}\quad = }&{\scriptstyle (\varphi\circ\varphi)\circ(\varphi\circ\varphi)=\varphi\{\varphi\}\{\varphi\{\varphi\}\}
=\varphi\{\varphi\{\varphi\{\varphi\}\}\}+\varphi\{\varphi,\varphi\{\varphi\}\}+\varphi\{\varphi\{\varphi\},\varphi\}. }
\end{array}
\eqno(\numb\label{eq103})
$$

In the computations we used the brace algebra identities~\eqref{eq89}.

The results are the same modulo 2.

\begin{figure}[!ht]
\unitlength=0.3em
\begin{center}
\begin{picture}(98,35)(-5,15)
  \put(1,45){
     \begin{picture}(0,0)
         \multiput(0,0)(0,-10){4}{\circle{4}}
         \multiput(0,-2)(0,-10){3}{\line(0,-1){6}}
          \put(-0.8,-1.2){$1$}
          \put(-0.8,-11.2){$2$}
          \put(-0.8,-21.2){$3$}
          \put(-0.8,-31.2){$4$}   
     \end{picture}
     }
    
   \put(6,35){$+$}   
   \put(12,35){$3$}
   \put(16,45){
      \begin{picture}(0,0)
          \put(5,0){\circle{4}}
          \put(0,-10){\circle{4}}
          \put(10,-10){\circle{4}}
          \put(10,-20){\circle{4}}
          \put(4.2,-1.2){$1$}
          \put(-0.8,-11.2){$2$}
          \put(9.2,-11.2){$3$}
          \put(9.2,-21.2){$4$}
          \put(1,-8.27){\line(1,2){3.2}}
          \put(9,-8.27){\line(-1,2){3.2}}
          \put(10,-12){\line(0,-1){6}}  
      \end{picture}
    }
    
   \put(32,35){$+$}   
   \put(38,35){$3$}
   \put(42,45){
      \begin{picture}(0,0)
          \put(5,0){\circle{4}}
          \put(0,-10){\circle{4}}
          \put(10,-10){\circle{4}}
          \put(0,-20){\circle{4}}
          \put(4.2,-1.2){$1$}
          \put(-0.8,-11.2){$2$}
          \put(-0.8,-21.2){$3$} 
          \put(9.2,-11.2){$4$}
          \put(1,-8.27){\line(1,2){3.2}}
          \put(9,-8.27){\line(-1,2){3.2}}
          \put(0,-12){\line(0,-1){6}}  
      \end{picture}
    }
 
   \put(58,35){$+$}   
   \put(64,35){$6$}
   \put(68,45){
      \begin{picture}(0,0)
         \put(10,0){\circle{4}}    
          \put(0,-10){\circle{4}}
          \put(10,-10){\circle{4}}  
          \put(20,-10){\circle{4}}
          \put(9.2,-1.2){$1$}
          \put(-0.8,-11.2){$2$}
          \put(9.2,-11.2){$3$}
          \put(19.2,-11.2){$4$}
          \put(1.42,-8.58){\line(1,1){7.16}}
          \put(18.58,-8.58){\line(-1,1){7.16}}
          \put(10,-2){\line(0,-1){6}}  
      \end{picture}
    }

\end{picture}
\end{center}
\caption{\footnotesize$\varphi^{[4]}=\varphi\{\varphi\{\varphi\{\varphi\}\}\}+3\varphi\{\varphi,\varphi\{\varphi\}\}+
3\varphi\{\varphi\{\varphi\},\varphi\}+6\varphi\{\varphi,\varphi,\varphi\}$}\label{fig9}
\end{figure}
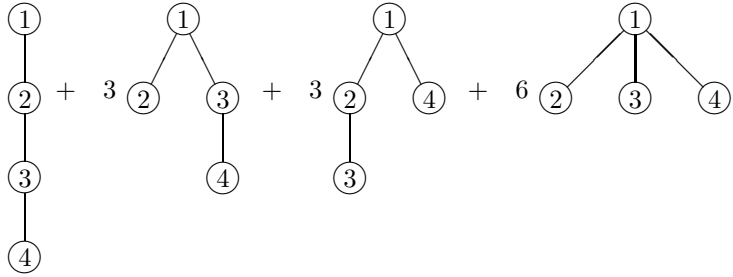

To visualize better the result of~\eqref{eq102} we will draw it as a sum of trees, see Figure~\ref{fig9}.

The numbers 1,2,3 and 4 are assigned in the order of appearance of $\varphi$ in the element of our sum. For
instance, for the element $\varphi\{\varphi,\varphi\{\varphi\}\}$: 1 is assigned to 
${\underline \varphi}\{\varphi,\varphi\{\varphi\}\}$, 2 is assigned to $\varphi\{{\underline \varphi},\varphi\{\varphi\}\}$, 3 is
assigned to $\varphi\{\varphi,{\underline \varphi}\{\varphi\}\}$, 4 is assigned to
$\varphi\{\varphi,\varphi\{{\underline\varphi}\}\}$. This enumeration corresponds to the standard recursion algorithm of
enumeration of vertices in a tree.

\medskip

Let us consider $\varphi^{[n]}$ for arbitrary $n$, $deg(\varphi)$ being odd. By means of the formula~\eqref{eq89}
we decompose $\varphi^{[n]}$ as a sum of similar trees with some coefficients.  By {\it tree} we mean a rooted
tree. The set of child vertices for each vertex is supposed to be ordered. All the trees will have $n$ vertices.
Consider any such tree $T$ and find its coefficient $O(T)$ in the decomposition. Let us label its vertices from
$1$ to $n$ as it was described before correspondingly to the standard recursion algorithm of enumeration of vertices. 
The structure of a rooted tree induces a partial order on the set of
vertices: the minimal element is the root that is always $1$; one says $i\prec j$ if the path from $1$ to $j$ passes
through $i$, see Definition~\ref{partialOrder}. It can be easily seen that the coefficient $O(T)$ is equal to the number of complete orders that
respect the partial order $\prec$.

The number $O(T)$ can be easily found. Denote by $ch_1(i),ch_2(i),\dots,ch_{k_i}(i)$, $i=1\dots n$, the child
vertices of vertex $i$; by $\#i$ --- the number of vertices that are greater or equal to vertex $i$ with respect to the partial order $\preceq$, we call this number 
$\#i$ {\it cardinality} of vertex $i$. Also we set
$$
\{  n_1,n_2,\dots,n_k\}:=\frac{(n_1+n_2+\dots+n_k)!}{n_1!n_2!\dots n_k!}. \eqno(\numb\label{eq104})
$$
It is the number of shuffles of $k$ sets of cardinalities $n_1$, $n_2$,$\dots$, $n_k$.

This expression is symmetric and  satisfies the following identity:
$$
\{  n_1,\dots,n_i,\dots,n_k\}\cdot\{  n_{i1},\dots,n_{im}\}=\{ 
n_1,\dots,n_{i-1},n_{i1},n_{i2},\dots,n_{im},n_{i+1},\dots,n_k\}, \eqno(\numb\label{eq104'})
$$
where $n_i=n_{i1}+n_{i2}+\dots+n_{im}$.

\medskip

One has
$$
O(T)=\prod_{i=1}^nO_i(T)=\prod_{i=1}^n\{  \#ch_1(i),\#ch_2(i),\dots,\#ch_{k_i}(i)\}.
\eqno(\numb\label{eq105})
$$
If the set of child vertices is empty for some vertex $i$, one sets $O_i(T):=1$.
\medskip

Let us prove the assertion  of the lemma.

Denote by $D_p(N)$ the sum of digits of the representation of $N$ in the base $p$ number system. It is easy to see
that $\{  n_1,n_2,\dots,n_k\}\not\equiv 0\mod p$ if and only if $D_p(n_1+n_2+\dots
+n_p)=D_p(n_1)+D_p(n_2)+\dots +D_p(n_k)$.

Now, consider $\varphi^{[pn]}$. Consider any tree $T$ with $pn$ vertices. We want to prove that $O(T)$ is
equal modulo $p$ to the coefficient of this tree in the decomposition of $(\varphi^{[p]})^{[n]}$. By the previous remark
$O(T)$ is not zero modulo $p$ if and only if

\medskip

($*$) For any vertex $i=1,2,..,pn$ we always have $D_p(\#ch_1(i)+\dots +\#ch_{k_i}(i))=D_p(\#ch_1(i))+\dots
+D_p(\#ch_{k_i}(i))$.

\medskip
In our proof we will need to examine only the last digit.

\medskip

\begin{figure}[!ht]
\begin{center}
\unitlength=0.3em
\begin{picture}(100,20)
   \put(0,20){
   \begin{picture}(0,20)
      \multiput(10,0)(-5,-10){3}{\circle{4}}
      \multiput(15,-10)(-5,-10){2}{\circle{4}}
      \put(20,-20){\circle{4}}     
      \put(9.2,-1.2){$1$}
      \put(4.2,-11.2){$2$}
      \put(-0.8,-21.2){$3$}
      \put(14.2,-11.2){$5$}
      \put(9.2,-21.2){$4$}
      \put(19.2,-21.2){$6$}
      \multiput(6,-8.27)(-5,-10){2}{\line(1,2){3.2}}
      \multiput(14,-8.27)(5,-10){2}{\line(-1,2){3.2}}
      \put(9,-18.27){\line(-1,2){3.2}}          
   \end{picture}
   }
   \put(26,10){\vector(1,0){8}}
   
   \put(37,20){
   \begin{picture}(0,20)
      \multiput(10,0)(-5,-10){3}{\circle{4}}
      \multiput(15,-10)(-5,-10){2}{\circle{4}}
      \put(20,-20){\circle{4}}     
      \put(9.2,-1.2){$1$}
      \put(4.2,-11.2){$2$}
      \put(-0.8,-21.2){$2$}
      \put(14.2,-11.2){$1$}
      \put(9.2,-21.2){$2$}
      \put(19.2,-21.2){$1$}
      \put(1,-18.27){\line(1,2){3.2}}
      \multiput(14,-8.27)(5,-10){2}{\line(-1,2){3.2}}
 
     \qbezier[20](6,-8.27)(7.5,-5) (8.9,-1.72)
      \put(9,-18.27){\line(-1,2){3.2}}          
   \end{picture}
   }
   
  \put(63,10){\vector(1,0){8}}
  
  \put(76,20){
  \begin{picture}(0,20)
    \put(2,-1){\circle{6}}
    \put(2,-18){\circle{6}}
    \thicklines
    \qbezier[20](1.9,-4)(1.9,-10)(1.9,-15)
    \put(0.9,-2.7){\Large $1$}
    \put(0.9,-19.7){\Large $2$}
  \end{picture}
  }

\end{picture}  
\end{center}
\caption{Assignment of indexes: $p=3$, $n=2$. Subtrees and quotient-tree.}\label{fig10}
\end{figure}
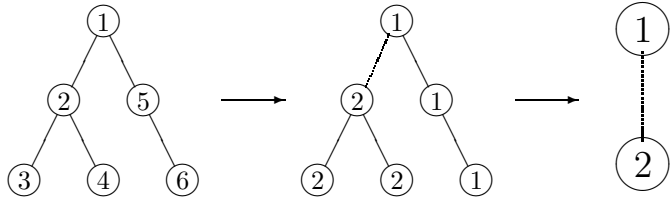

Remind that $pn$ vertices of our tree $T$ are enumerated from $1$ to $pn$. We will also index them by numbers from $1$ to
$n$ as follows. We pass all the vertices accordingly to their enumeration. We assign 1 to the root. If for vertex
$i$ we have cardinality $\#i$ is not a number divisible by $p$, we assign to this vertex the same index as for
its parent vertex, if $\#i$ is divisible by $p$ we assign the next number that was not yet used. In the case
tree $T$ satisfies the condition ($*$) we will use exactly $n$ numbers (otherwise this number can be less), see
Figure~\ref{fig10}.

If $O(T)\not\equiv 0\mod p$, then for any number $j=1\dots n$ the vertices with the same index $j$ form a
subtree with $p$ vertices. Denote this subtree by $T_j$. We draw the edges of trees $T_j$ by continuous
lines, all the other edges --- by dotted lines. Consider quotient-tree $H$ obtained from $T$ by contraction of
continuous edges, see Figure~\ref{fig10}. Tree $H$ has $n$ vertices. We will prove that
$$
O(T)\equiv O(H)\cdot O(T_1)\cdot O(T_2)\cdot\ldots\cdot O(T_n)\mod p. \eqno(\numb\label{eq106})
$$
Evidently, this will prove the lemma: tree $H$ and coefficient $O(H)$ correspond to tree decomposition of
$[n]$-operation; trees $T_1,\dots,T_n$ and coefficients $O(T_1),\dots ,O(T_n)$ correspond to tree decomposition of
$[p]$-operation.

Consider any vertex $i=1\dots pn$ indexed by some $j=1\dots n$, suppose its number in the tree $T_j$ is $s$.
It has child vertices of two types: those that are  indexed by $j$, suppose their cardinalities are
$a_1p+b_1,a_2p+b_2,\dots,a_\ell p+b_\ell$; and those that have another index, suppose their cardinalities
are $a_{\ell+1}p,a_{\ell+2}p,\dots,a_Lp$. The numbers $a_1,\dots,a_L,b_1,\dots,b_l$ are positive and satisfy
$A:=a_1+a_2+\dots +a_L<n$, $B:=b_1+b_2+\dots+b_l<p$.

Note that $\{  mp,b\}\equiv 1\mod p$, if $b<p$.

By the previous remark and by propriety~\eqref{eq104'} we get:
$$
\begin{array}{rcl}
O_i(T)&=&\{  a_1p+b_1,a_2p+b_2,\dots,a_lp+b_l,a_{l+1}p,\dots,a_Lp\}\\
&\equiv&\{ a_1p,b_1,a_2p,b_2,\dots,a_lp,b_l,a_{l+1}p,a_{l+2}p,\dots,a_Lp\}\\
&=&\{a_1p,a_2p,\dots,a_Lp,b_1,b_2,\dots,b_l\}\\
&=&\{  Ap,B\}\cdot \{  a_1p,\dots a_Lp\}\cdot\{  b_1,\dots,b_l\}\\
&\equiv&\{  a_1p,\dots a_Lp\}\cdot\{  b_1,\dots,b_l\}\\
&\equiv&\{  a_1,\dots a_L\}\cdot\{  b_1,\dots,b_l\}\\
&=&\{  a_1,\dots a_L\}\cdot O_s(T_j)\mod p.
\end{array}
$$
The next to last equality was obtained due to the fact that $\{  a_1p,\dots a_Lp\}\equiv\{ 
a_1,\dots a_L\}\mod p$.

If we take the product of expressions $\{  a_1,\dots a_L\}$ for all vertices of the tree $T_j$ then we
obtain exactly $O_j(H)$ (we will need to apply many times the identity~\eqref{eq104'}). This proves~\eqref{eq106}
and thus the lemma. \nbox

\bigskip

\noindent {\bf Proof of Lemma~\ref{l1010}:}. First we need to understand how $\varphi^{[n]}$ is decomposed as a
sum of trees if $deg(\varphi)$ is even. More precisely we need to find the appropiate coefficients $O^-(T)$ of all trees $T$ with
$n$ vertices in this case. Combinatorically we will have the same number $O(T)$ of summands that give the same tree $T$ (the number
of complete orders that respect the partial order of vertices induced by the structure of a tree). But in this
situation not all of the summands are "1" ---  they are $\pm 1$. The sign is the sign of the  permutation  corresponding to
a complete order (remind that we have the fixed enumeration, with respect to this enumeration one defines
permutations).

Denote by
$$
\{  n_1,n_2,\dots,n_k\}_{-1} \eqno(\numb\label{eq104''})
$$
the superanalog of the expression~\eqref{eq104}, see Appendix~\ref{C}. This number is the number of shuffles of sets of
cardinalities $n_1,\dots,n_k$, that provide even permutations of elements, minus the number of shuffles of the
same sets providing odd permutations. We will call it the {\it number of shuffles of sets with anticommuting
elements}. This number is
$$
\{  n_1,n_2,\dots,n_k\}_{-1}=
\begin{cases}
0, \, \text{\footnotesize if at least two numbers $n_i$ and $n_j$, $\{i,j\}\subset\{1,\dots,k\}$, are odd;}\\
\{  \left[\frac{n_1}2\right],\left[\frac{n_2}2\right],\dots,\left[\frac{n_k}2\right]\},\, \text{\footnotesize otherwise,}
\end{cases}
\eqno(\numb\label{eq107})
$$
see Appendix~\ref{C}. We see from this formula that expression~\eqref{eq104''} is symmetric: it does not depend on the
order of numbers $n_1,\dots,n_k$. 

\smallskip

Analogously to~\eqref{eq105} we have:
$$
O^-(T)=\prod_{i=1}^nO^-_i(T)=\prod_{i=1}^n\{  \#ch_1(i),\#ch_2(i),\dots,\#ch_{k_i}(i)\}_{-1}.
\eqno(\numb\label{eq109})
$$

\bigskip

We are ready to prove the lemma.  Consider any tree $T$ with $2n$ vertices. By~\eqref{eq107} inequality
$O^-(T)\neq 0$ implies

\medskip

($**$) for any vertex $i=1,\dots,2n$ all the cardinalities $\#ch_1(i),\#ch_2(i),\dots,\#ch_{k_i}(i)$ (of its child
vertices) except probably one are even.

\medskip

Suppose our tree $T$ satisfies this condition ($**$). The vertices of $T$ are enumerated from $1$ to $2n$. By means of the 
algorithm descibed in the proof of the previous lemma we will
index them by numbers $1,2,\dots,n$:  We assign 1 to the root. If for vertex
$i$ we have cardinality $\#i$ is odd, we assign to this vertex the same index as for its parent vertex, if $\#i$
is even we assign the next number that was not yet used. In the case tree $T$ satisfies the condition ($**$) we
will use exactly $n$ numbers (otherwise this number will be less), see Figure~\ref{fig11}. Finally we obtain
subtrees $T_1,\dots,T_{n}$ with 2 vertices --- they are all 
\unitlength=0.08em
\begin{picture}(10,20)(0,3)
  \multiput(5,3)(0,14){2}{\circle{6}}
  \put(5,6){\line(0,1){8}}
\end{picture}, 
and a quotient tree $H$ with $n$ vertices, see
Figure~\ref{fig11}.

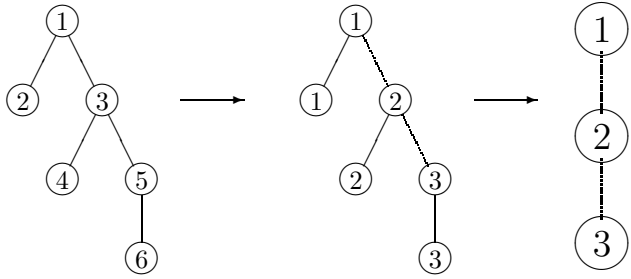
\begin{figure}[!ht]
\begin{center}
\unitlength=0.3em
\begin{picture}(100,30)
   \put(0,30){
   \begin{picture}(0,30)
      \multiput(10,0)(-5,-10){2}{\circle{4}}
      \multiput(15,-10)(-5,-10){2}{\circle{4}}
      \put(20,-20){\circle{4}}
      \put(20,-30){\circle{4}}     
      \put(9.2,-1.2){$1$}
      \put(4.2,-11.2){$2$}
      \put(14.2,-11.2){$3$}
      \put(9.2,-21.2){$4$}
      \put(19.2,-21.2){$5$}
      \put(19.2,-31.2){$6$}
      \multiput(6,-8.27)(5,-10){2}{\line(1,2){3.2}}
      \multiput(14,-8.27)(5,-10){2}{\line(-1,2){3.2}}
      \put(20,-22){\line(0,-1){6}}          
   \end{picture}
   }
   \put(26,20){\vector(1,0){8}}
   
   \put(37,30){
    \begin{picture}(0,30)
      \multiput(10,0)(-5,-10){2}{\circle{4}}
      \multiput(15,-10)(-5,-10){2}{\circle{4}}
      \put(20,-20){\circle{4}}
      \put(20,-30){\circle{4}}     
      \put(9.2,-1.2){$1$}
      \put(4.2,-11.2){$1$}
      \put(14.2,-11.2){$2$}
      \put(9.2,-21.2){$2$}
      \put(19.2,-21.2){$3$}
      \put(19.2,-31.2){$3$}
      \multiput(6,-8.27)(5,-10){2}{\line(1,2){3.2}}
      \multiput(14,-8.27)(5,-10){2}{\qbezier[20](0,0)(-1.6,3.2)(-3.2,6.4)}
      \put(20,-22){\line(0,-1){6}}          
   \end{picture}
  
   }
   
  \put(63,20){\vector(1,0){8}}
  
  \put(76,30){
  \begin{picture}(0,30)
    \put(2,-1){\circle{6}}
    \put(2,-14.5){\circle{6}}
    \put(2,-28){\circle{6}}
    \thicklines
    \multiput(1.9,-4)(0,-13.5){2}{\qbezier[15](0,0)(0,-4)(0,-7.5)}

    \put(0.9,-2.7){\Large $1$}
    \put(0.9,-16.2){\Large $2$}
    \put(0.9,-29.7){\Large $3$}
  \end{picture}
  }

\end{picture}  
\end{center}
\caption{Assignment of indexes. Subtrees and quotient-tree.}\label{fig11}
\end{figure}

Consider any subtree $T_j$, $j=1,\dots,n$. It has two vertices: the upper one $i_1$ and its child $i_2$, where
$i_1,i_2\in\{1,\dots,2n\}$. Let us compute $O^-_{i_1}(T)$, $O^-_{i_2}(T)$. Suppose cardinalities of the child
vertices of $i_1$ are $2a_1+1,2a_2,2a_3,2a_4,\dots,2a_\ell$ (the only odd cardinality is that of $i_2$).
Respectively for $i_2$ this cardinalities are $2a_{11},2a_{12},\dots,2a_{1t}$ (where
$2a_{11}+2a_{12}+\dots+2a_{1t}=2a_1$). Thus
$$
\begin{array}{rl}
O^-_{i_1}(T)\cdot O^-_{i_2}(T)&=\{  2a_1+1,2a_2,2a_3,2a_4,\dots,2a_\ell\}_{-1}\cdot\{ 
2a_{11},2a_{12},\dots,2a_{1t}\}_{-1}\\
&=\{  a_1,a_2,a_3,\dots,a_\ell\}\cdot\{ 
a_{11},a_{12},\dots,a_{1t}\}\\
&=\{  a_{11},a_{12},\dots,a_{1t},a_2,a_3,\dots,a_\ell\}\\
&=O_j(H)=O_j(H)\cdot O^-(T_j).
\end{array}
$$
The second equality follows from~\eqref{eq107}. The last equality is due to $O^-(\unitlength=0.08em
\begin{picture}(10,20)(0,5)
  \multiput(5,3)(0,14){2}{\circle{6}}
  \put(5,6){\line(0,1){8}}
\end{picture})=1$.

As a consequence we obtain  
$$
O^-(T)=O(H)=O(H)\cdot O^-(T_1)\cdot O^-(T_2)\cdot\ldots \cdot O^-(T_n).
$$
This proves
the lemma --- tree $H$ and coefficient $O(H)$ correspond to $[n]$-operation, trees $T_1,\dots,T_n$ and
coefficients $O^-(T_1),\dots,O^-(T_n)$ correspond to $[2]$-operation. \nbox

\section{Pairing}\label{s11}
Any diagram of bigrading~$(i,i+1)$ is a cycle in complex $CTD$ homologous to  diagram~\eqref{eq72} taken with
some coefficient, see Theorem~\ref{t71}. In this section we find out this coefficient.

Consider any alternated $T$-diagram $A$ of complexity $(i,i+1)$. It has only one minimal component, \ie  only one
alternated tree. By abuse of the language we designate this alternated tree also by $A$. Let us assign to $A$ a
tree in the sense of previous section. This tree $T(A)$ will have $i$ vertices. We will do it recursively. If
$i=1$ there is only one diagram  
\unitlength=0.15em
\begin{picture}(20,10)
     \put(0,0){\line(1,0){20}}
      \qbezier(5,0)(10,9)(15,0)
\end{picture}. 
    We assign to it one vertex tree: 
    $T( \begin{picture}(20,10)
     \put(0,0){\line(1,0){20}}
      \qbezier(5,0)(10,9)(15,0)
    \end{picture})=
    \begin{picture}(10,10) \put(4,2.4){\circle{4}} \end{picture}$. Let $i>1$. Since $A$ is
alternated it always contains the edge $(1,i+1)$ joining the extremal vertices. If we remove this edge we obtain
whether two disconnected non-trivial alternated trees $A_1$ and $A_2$, whether an alternated tree $A_1$ and a single point. In
the first situation we set
$$
T(A)=
\unitlength=0.15em
\begin{picture}(45,15)
 \put(19,10){\circle{4}}
 \put(18,8.3){\line(-4,-3){10}}
 \put(20,8.3){\line(4,-3){10}}
 \put(0,-5){$T(A_{1})$}
 \put(25,-5){$T(A_{2})$}
\end{picture}.
$$
In the second situation we set
$$
T(A)=
\unitlength=0.15em
\begin{picture}(23,15)
 \put(7.5,10){\circle{4}}
 \put(7.5,8){\line(0,-1){7}}
 \put(1.5,-5){$T(A_{1})$}

\end{picture}.
$$

\medskip

\noindent For instance, in the following situation this algorithm gives:
$$
T(
\unitlength=0.17em
\begin{picture}(62,35)(-6,3)
 \put(-3,0){\line(1,0){56}}
 \qbezier(0,0)(25,30)(50,0)
 \qbezier(0,0)(20,20)(40,0)
 \qbezier(0,0)(15,15)(30,0)
 \qbezier(0,0)(10,10)(20,0)
 \qbezier(10,0)(20,10)(30,0)

\end{picture})\quad =
\begin{picture}(40,30)(-10,10)
 \put(10,30){\circle{5}}
 \put(10,20){\circle{5}}
 \put(10,10){\circle{5}}
 \put(0,0){\circle{5}}
 \put(20,0){\circle{5}}
 \multiput(10,27.5)(0,-10){2}{\line(0,-1){5}}
 \put(8.23,8.23){\line(-1,-1){6.5}}
 \put(11.67,8.23){\line(1,-1){6.5}}
\end{picture}
$$

\vspace{4mm}

\begin{proposition}\label{p111}
Any alternated $T$-diagram $A$ of bigrading $(i,i+1)$ is homologous to zero if it has intersecting edges. If $A$
has no intersecting edges, it is homologous to diagram~\eqref{eq72}

(i) in complex $CTD^{even}$: with coefficient $\pm O(T(A))$;

(ii) in complex $CTD^{odd}$: with coefficient $\pm O^-(T(A))$.

\smallskip

The sign depends on the orientation of $A$. The coefficient must be taken modulo the order of the corresponding
homology group.\nbox
\end{proposition}

Precision: we say a diagram has intersecting edges if it has vertices $a<b<c<d$ on the line $\R^1$, and has edges
$(a,c)$, $(b,d)$. For the definition of numbers $O(T)$, $O^-(T)$ see the previous section.

\medskip

\noindent{\bf Proof of Proposition~\ref{p111}:} Immediately follows from the duality of the basis of alternated trees
to that of monotone brackets and also from the considerations of the previous section. \nbox

\section{How the homology bialgebra of $DHATD$ is related to that of $DHAT_0D$}\label{s14}

In this section we will completely describe the relation between the homology bialgebra of $DHATD$ and that of
$DHAT_0D$ for any ring $\kk$ of coefficient. The results of this section (namely Theorem~\ref{t116}) are given
without proof. The proof will be given elsewhere.

First we need to fix some new objects and notations.

The following definition is a generalization of~\ref{d75}.

\begin{definition}\label{d111}
{\rm A {\it divided product} $\langle  D_1,D_2,\dots,D_n\rangle $ of $T_{*}$/$T$/$T_0$-diagrams $D_1,D_2,\dots,D_n$ is the sum of those
elements in the shuffle product $D_1*D_2*\ldots *D_n$ that have the left-most point of $D_i$ on the left from the
left-most point of $D_{i+1}$ for all $i=1,\dots,n-1$.} \nbox
\end{definition}
We extend these operations as multilinear operations on the space of $T_{*}$-diagrams (resp. $T$-diagrams and $T_0$-diagrams).

We will denote by
$$
\langle\,  \rangle 
$$
the trivial diagram --- the unity of algebras $DHAT_{*}D$, $DHATD$ and $DHAT_0D$.

\begin{lemma}\label{l112}
{\footnotesize
\begin{gather*}
\Delta\langle D_1,D_2,\dots,D_n\rangle =\sum_{i=0}^{n} \langle D_1,D_2,\dots,D_i\rangle \otimes\langle D_{i+1},D_{i+2},\dots,D_n\rangle .\\
\langle D_1,D_2,\dots,D_n\rangle * \langle D_{n+1},D_{n+2},\dots,D_{n+m}\rangle =\sum_{\sigma\in S(n,m)}
(-1)^{s(\sigma)}\langle D_{\sigma(1)},D_{\sigma(2)},\dots,D_{\sigma(n+m)}\rangle ,
\end{gather*}
}
where $S(n,m)$ is a subset in the symmetric group $S_{n+m}$ whose elements are shuffles of  $1,2,\dots,n$ with
$n+1,n+2,\dots,n+m$; $(-1)^{s(\sigma)}$ is the sign of the induced permutation of odd elements among
$D_1,D_2,\dots,D_{n+m}$. \nbox
\end{lemma}

Proof of this lemma is a direct check.

\medskip

Consider subspaces $\ZZ^{even}\subset CTD^{even}$, $\ZZ^{odd}\subset CTD^{odd}$ that are spanned by the elements
$$
\langle Z_{k_1},Z_{k_2},\dots,Z_{k_l}\rangle .
$$
Remind, $Z_i$ is the diagram~\eqref{eq72}. By the previous lemma these subspaces are Hopf subalgebras of
$DHATD^{even}$, resp. $DHATD^{odd}$. The following lemma generalizes~\eqref{eq77},~\eqref{eq78} and shows that
$\ZZ^{even}$, $\ZZ^{odd}$
 are differential Hopf subalgebras:

\begin{lemma}\label{l113}
For even $d$ one has
$$
\begin{array}{l}
\partial\langle Z_{k_1},Z_{k_2},\dots,Z_{k_l}\rangle =\\
\sum_{i=1}^{l-1}(-1)^{i-1}{{k_i+k_{i+1}}\choose{k_i}}
\langle Z_{k_1},Z_{k_2},\dots,Z_{k_{i-1}},Z_{k_i+k_{i+1}},Z_{k_{i+2}},\dots,Z_{k_l}\rangle .
\end{array}
\eqno(\numb)\label{eq114}
$$
For odd $d$ one has
$$
\begin{array}{l}
\partial\langle Z_{k_1},Z_{k_2},\dots,Z_{k_l}\rangle =\\
\sum_{i=1}^{l-1}(-1)^{i-1+k_1+k_2+\dots+k_i}{{k_i+k_{i+1}}\choose{k_i}}_{-1}
\langle Z_{k_1},Z_{k_2},\dots,Z_{k_{i-1}},Z_{k_i+k_{i+1}},Z_{k_{i+2}},\dots,Z_{k_l}\rangle .
\end{array}
\eqno(\numb)\label{eq115}
$$
$\Box$
\end{lemma}

Relation between $DHATD$ and $DHAT_0D$ is given by the following theorem.

\begin{theorem}\label{t116}
The morphisms
\begin{gather}
\mu:\ZZ^{even}\otimes DHAT_0D^{even}\rightarrow DHATD^{even},\\
\mu:\ZZ^{odd}\otimes DHAT_0D^{odd}\rightarrow DHATD^{odd}
\end{gather}
are quasi-isomorphisms of differential Hopf algebras ($\mu$ is the shuffle multiplication in $DHATD$). \nbox
\end{theorem}
Without proof. The proof of this theorem is more or less geometrical and is absolutely different from the considerations of this article. This
theorem together with Theorem~\ref{t71} provide another proof of Theorem~\ref{t72}.

\smallskip

For instance if our ring of coefficients $\kk$ is a field, then this theorem says that the homology bialgebra of
$DHATD$ is a tensor product of homology bialgebra of $DHAT_0D$ to that of $\ZZ$.

Complexes $\ZZ^{even}$, $\ZZ^{odd}$  are well known, see~\cite{F,Vain,V0,Ma}. In Appendix~D we give a summary
describing the structure of their homology bialgebras.

For a fixed complexity $i$, complex $\ZZ^{odd}$ (resp. $\ZZ^{even}$) computes  cohomology groups $H^*(Br(i),\Z)$
of the braid group with $i$ strings (resp. cohomology groups $H^*(Br(i),\pm\Z)$ of this braid group with
coefficients in its sign representation $\pm\Z$). Indeed, classifying space of the braid group $Br(i)$ is the
configuration  space $B(\C^1,i)$ (space of cardinality $i$ subsets of $\C^1$). By Poincar\'e duality one has
\begin{gather*}
H^*(B(\C^1,i),\Z)\simeq \bar H_{2i-*}(B(\C^1,i),\Z),\\
H^*(B(\C^1,i),\pm\Z)\simeq \bar H_{2i-*}(B(\C^1,i),\pm\Z).
\end{gather*}
($\bar H_*(\, .\, ,L)$ designates locally finite homology groups with coefficients in a local system $L$.) One
point compactification $\overline{B(\C,i)}$ of the configuration space has a natural cell decomposition that is
defined as follows. Let  $\xi=\{z_1,z_2,\dots,z_i\}$ be a point of $B(\C^1,i)$. We will assign to $\xi$ its {\it
index} --- system of numbers $(k_1,k_2,\dots,k_l)$ satisfying $k_1+k_2+\dots +k_l=i$, where $k_1$ is the number of
elements of $\xi$ with the minimal value of the  real component $\Re (z)$; $k_2$ is the number of elements of
$\xi$ with next value of $\Re (z)$, and so on$\dots$ Points with the same index $(k_1,k_2,\dots,k_l)$ form a cell,
that we denote by $e(k_1,k_2,\dots,k_l)$. 

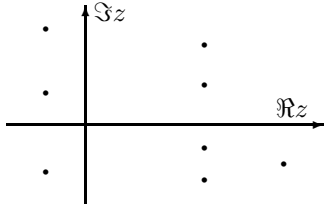
\begin{figure}[h!t]
\unitlength=0.3em
\begin{center}
\begin{picture}(40,25)

\put(0,10){\vector(1,0){40}}
\put(10,0){\vector(0,1){25}}
\put(5,4){\circle*{0.5}}
\put(5,14){\circle*{0.5}}
\put(5,22){\circle*{0.5}}

\put(25,15){\circle*{0.5}}
\put(25,20){\circle*{0.5}}
\put(25,7){\circle*{0.5}}
\put(25,3){\circle*{0.5}}

\put(35,5){\circle*{0.5}}
\put(34,11){$\Re z$}
\put(11,23){$\Im z$}
\end{picture}
\end{center} 
\caption{Point of the cell $e(3,4,1)$ of $\overline{B(\C,8)}$.}  
\end{figure}

All such cells together with the infinite point provide a cell
decomposition of $\overline{B(\C^1,i)}$. These cells bound to each other exactly by the 
rule~\eqref{eq115}, see~\cite{Vain}
($\langle Z_{k_1},\dots,Z_{k_l}\rangle $  should be replaced by $e(k_1,\dots,k_l)$).

The differential of the cells in the local system $\pm\Z$ is described by the rule~\eqref{eq114}, see~\cite{V0,Ma}.

\medskip

If we consider complexes $\ZZ^{even}$, $\ZZ^{odd}$ as differential Hopf algebras and fix as only grading
$deg=p+q=(d-1)i-j$, then each complex computes the cohomology bialgebra
$H^*(\Omega^2(S^{d-1}),\Z)$ of the double loop space of $(d-1)$-dimensional sphere. To see this one can use
Vassiliev's approach of discriminants, see~\cite{V4}. Connection with the homology of braid groups is a particular case
of the well known Snaith decomposition formula, see~\cite{Snaith}:
$$
H^t(\Omega^mS^n)\simeq\oplus_{i=1}^\infty H^{t-i(n-m)}(B(\R^m,i),(\pm\Z)^{\otimes(n-m)})
$$
for $m=2$, $n=d-1$. ($B(\R^n,i)$ denotes the space of cardinality $i$ subsets of $\R^n$.)

\medskip

Theorem~\ref{t116} has a very simple geometrical meaning. The space $Emb=\K\setminus\Sigma$ of long 
knots is a subspace of the space
$Imm$ of long immersions --- immersions with a fixed behavior at infinity:
$$
Emb\hookrightarrow Imm.
$$
 The homotopy fiber $Emb^+$ of this
inclusion can be regarded as a Serre fibration over $Emb$:
$$
\Pi:Emb^+\to Emb.
$$
Preimage $\Pi^{-1}(f)$ of any knot $f\in Emb$ is the space $\Omega (Imm;f,l)$ of pathes in $Imm$ that start in $f$
and end in the fixed linear embedding $l$. Obviously, $\Pi^{-1}(f)$ is homotopy equivalent to the loop space
$\Omega(Imm)$. On the other hand, the space $Imm$ is homotopy equivalent to
$\Omega(\R^d\setminus\{0\})\simeq\Omega S^{d-1}$ (we consider the value of derivative $f'(t)$, $f\in Imm$, to
obtain a map in one direction).  Thus we have $\Pi^{-1}(k)\simeq \Omega^2(S^{d-1})$.
Actually, complex $DHATD$ (resp. $DHAT_0D$) computes the second term of  D.~Sinha's spectral sequence converging to the 
(co)homology of
the space $Emb^+$ (resp. $Emb$), \cf~\cite{Vol,Sinha2}.  A crucial point is that $\Pi$ is a trivial fiber bundle, \ie 
$Emb^+\simeq Emb\times \Omega^2 S^{d-1} $,  since 
the inclusion $Emb\hookrightarrow Imm$ is a contractible map, 
see~\cite{Sinha2}. Therefore, Theorem~\ref{t116} confirms the conjecture that Sinha's spectral sequence stabilizes at the second term.

\appendix

\part*{Appendixes}

In Appendix~\ref{A} we describe bases of complexes $CTD$, $(\Poiss^{Norm},\partial)$, $CT_{*}D$, $(\BV^{Norm},\partial)$, $CT_{0}D$,
$(\Poiss^{zero},\partial)$ for small complexities $i$.

Appendix~\ref{B} contains results of computer calculations of the homology of $CTD$ and $CT_{0}D$.

Appendix~\ref{C} gives shuffle combinatorial formulas (like binomial) for anticommuting elements.

In Appendix~\ref{D} we describe the  homology of differential bialgebra $\ZZ$. This bialgebra is well known and its homology bialgebra
is $H^*(\Omega^2S^{d-1})$, $d\ge 4$.

Appendix~\ref{E} describes homology bialgebra $H_{*}(\Omega^2S^{d-1})$ and provides explicit formulas for inclusion
of $H_{*}(\Omega^2S^{d-1})$ in Hochschild homology of Poisson or Gerstenhaber algebras operad.

\section{Complexes of $T_*/T/T_0$-diagrams for small complexities $i$}\label{A}
Complexes of $T_*/T/T_0$-diagrams are finite for each complexity $i$. The aim of this section is to describe the
bases of these complexes for small $i$. In particular, this demonstrates how complex $CT_0D$ simplifies computations of the
homology groups of $CT_*D$. We describe also the dual bases in the dual complexes.

\subsection{Alternated $T$-diagrams of complexities $i=1,2$}
In complexity $i=1$ one has only one alternated $T$-diagrams:
$$
\unitlength=0.3em
\begin{picture}(20,7)
     \put(0,0){\line(1,0){20}}
      \qbezier(5,0)(10,7)(15,0)
\end{picture}
$$
The dual element is $[x_{1},x_{2}]$.

\vspace{3mm}

In complexity $i=2$ one has five alternated $T$-diagrams:
$$
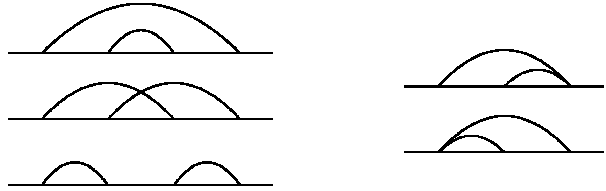
\unitlength=0.25em
\begin{picture}(90,30)
  \put(0,0){
  \begin{picture}(40,10)
    \put(0,0){\line(1,0){40}}
    \multiput(5,0)(20,0){2}{\qbezier(0,0)(5,7)(10,0)}
  \end{picture}
  }
  
  \put(0,10){
  \begin{picture}(40,10)
    \put(0,0){\line(1,0){40}}
    \multiput(5,0)(10,0){2}{\qbezier(0,0)(10,11)(20,0)}
  \end{picture}
  }
  
  \put(0,20){
  \begin{picture}(40,10)
    \put(0,0){\line(1,0){40}}
    \put(5,0){\qbezier(0,0)(15,15)(30,0)}
    \put(15,0){\qbezier(0,0)(5,7)(10,0)}
  \end{picture}
  }
  
  \put(60,5){
  \begin{picture}(30,10)
    \put(0,0){\line(1,0){30}}
    \put(5,0){\qbezier(0,0)(5,5)(10,0)}
    \put(5,0){\qbezier(0,0)(10,11)(20,0)}
  \end{picture}
  }

  \put(60,15){
  \begin{picture}(30,10)
    \put(0,0){\line(1,0){30}}
    \put(15,0){\qbezier(0,0)(5,5)(10,0)}
    \put(5,0){\qbezier(0,0)(10,11)(20,0)}
  \end{picture}
  }
  
\end{picture}
$$
The dual monotone bracket diagrams are $[x_{1},x_{4}]\cdot[x_{2},x_{3}]$, \,  $[x_{1},x_{3}]\cdot[x_{2},x_{4}]$, \, 
$[x_{1},x_{2}]\cdot[x_{3},x_{4}]$
 \, and \,
$[x_{1},[x_{2},x_{3}]]$, \, $[[x_{1},x_{2}],x_{3}]$ respectively.

\subsection{Alternated $T_*$-diagrams of complexities $i=1,2$}
In complexity $i=1$ one has two alternated $T_*$-diagrams:
$$
\unitlength=0.3em
\begin{picture}(20,10)
     \put(0,0){\line(1,0){20}}
      \qbezier(5,0)(10,8)(15,0)
\end{picture}
\qquad\qquad
\begin{picture}(15,10)
     \put(0,0){\line(1,0){15}}
      \put(6.5,-0.8){$*$}
\end{picture}
$$
The dual elemenents in the $\BV^{Norm}$ are $[x_{1},x_{2}]$ and $\delta x_{1}$.

In complexity $i=2$ one has eleven alternated $T_*$-diagrams:
$$
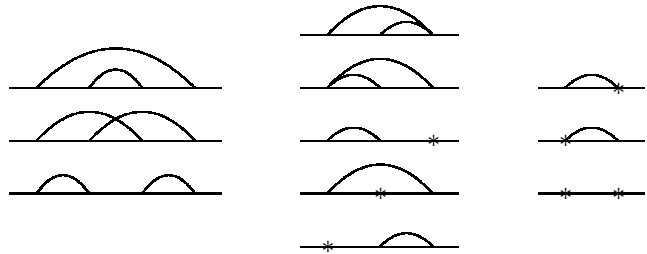
\unitlength=0.2em
\begin{picture}(120,50)
  \put(0,10){
  \begin{picture}(40,10)
    \put(0,0){\line(1,0){40}}
    \multiput(5,0)(20,0){2}{\qbezier(0,0)(5,7)(10,0)}
  \end{picture}
  }
  
  \put(0,20){
  \begin{picture}(40,10)
    \put(0,0){\line(1,0){40}}
    \multiput(5,0)(10,0){2}{\qbezier(0,0)(10,11)(20,0)}
  \end{picture}
  }
  
  \put(0,30){
  \begin{picture}(40,10)
    \put(0,0){\line(1,0){40}}
    \put(5,0){\qbezier(0,0)(15,15)(30,0)}
    \put(15,0){\qbezier(0,0)(5,7)(10,0)}
  \end{picture}
  }
  
  \put(55,30){
  \begin{picture}(30,10)
    \put(0,0){\line(1,0){30}}
    \put(5,0){\qbezier(0,0)(5,5)(10,0)}
    \put(5,0){\qbezier(0,0)(10,11)(20,0)}
  \end{picture}
  }

  \put(55,40){
  \begin{picture}(30,10)
    \put(0,0){\line(1,0){30}}
    \put(15,0){\qbezier(0,0)(5,5)(10,0)}
    \put(5,0){\qbezier(0,0)(10,11)(20,0)}
  \end{picture}
  }

  \put(55,20){
   \begin{picture}(30,10)
     \put(0,0){\line(1,0){30}}
      \qbezier(5,0)(10,5)(15,0)
      \put(24,-1.2){$*$}
   \end{picture}
   }

  \put(55,10){
   \begin{picture}(30,10)
     \put(0,0){\line(1,0){30}}
      \qbezier(5,0)(15,11)(25,0)
      \put(14,-1.2){$*$}
   \end{picture}
   }

  \put(55,0){
   \begin{picture}(30,10)
     \put(0,0){\line(1,0){30}}
      \qbezier(15,0)(20,5)(25,0)
      \put(4,-1.2){$*$}
   \end{picture}
   }

  \put(100,30){
  \begin{picture}(20,10)
    \put(0,0){\line(1,0){20}}
    \qbezier(5,0)(10,5)(15,0)
     \put(14,-1.2){$*$}
   \end{picture}  
   }
  
  \put(100,20){
  \begin{picture}(20,10)
    \put(0,0){\line(1,0){20}}
    \qbezier(5,0)(10,5)(15,0)
     \put(4,-1.2){$*$}
   \end{picture}  
   }
   
     \put(100,10){
  \begin{picture}(20,10)
    \put(0,0){\line(1,0){20}}
    \put(4,-1.2){$*$}
     \put(14,-1.2){$*$}
   \end{picture}  
   }

\end{picture}
$$

The dual elements in $\BV^{Norm}$ are:
\begin{center}
\begin{tabular}{ccccc}
&& $[x_{1},[x_{2},x_{3}]]$&& \\
$[x_{1},x_{4}]\cdot[x_{2},x_{3}]$ && $[[x_{1},x_{2}],x_{3}]$ && $[x_{1},\delta x_{2}]$ \\
$[x_{1},x_{3}]\cdot[x_{2},x_{4}]$ && $[x_{1},x_{2}]\cdot \delta x_{3}$ && $[\delta x_{1},x_{2}]$ \\
$[x_{1},x_{2}]\cdot[x_{3},x_{4}]$ && $[x_{1},x_{3}]\cdot \delta x_{2}$ && $\delta x_{1} \cdot \delta x_{2}$ \\
&& $\delta x_{1}\cdot  [x_{2},x_{3}]$&& \\
\end{tabular}
\end{center}

\subsection{Alternated $T_0$-diagrams of complexities $i=1,2,3$}
There is no alternated $T_0$-diagrams of complexity $i=1$.

In complexity $i=2$ one has only one alternated $T_0$-diagram:
$$
\unitlength=0.2em
\begin{picture}(40,10)
    \put(0,0){\line(1,0){40}}
    \multiput(5,0)(10,0){2}{\qbezier(0,0)(10,11)(20,0)}
  \end{picture}
$$
The dual bracket diagram is $[x_{1},x_{3}]\cdot [x_{2},x_{4}]$.

In complexity $i=3$ one has twelve  alternated $T_0$-diagrams:
$$
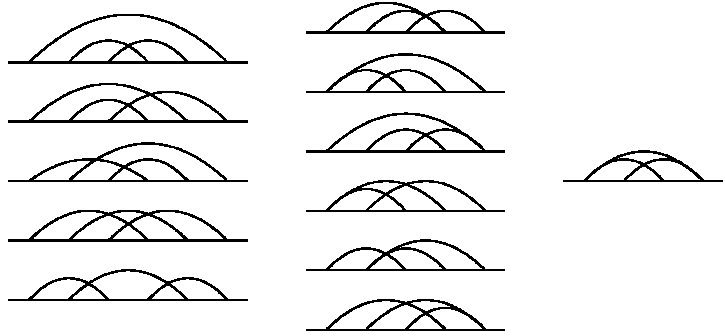
\unitlength=0.15em
\begin{picture}(180,90)
  \put(0,67.5){
  \begin{picture}(60,15)
    \put(0,0){\line(1,0){60}}
    \qbezier(5,0)(30,24)(55,0)
    \qbezier(15,0)(25,11)(35,0)
    \qbezier(25,0)(35,11)(45,0)
  \end{picture}
  }
  
  \put(0,52.5){
  \begin{picture}(60,15)
    \put(0,0){\line(1,0){60}}
    \qbezier(5,0)(25,19)(45,0)
    \qbezier(15,0)(25,11)(35,0)
    \qbezier(25,0)(40,15)(55,0)
  \end{picture}
  }
  
  \put(0,37.5){
  \begin{picture}(60,15)
    \put(0,0){\line(1,0){60}}
    \qbezier(5,0)(20,11)(35,0)
    \qbezier(15,0)(35,19)(55,0)
    \qbezier(25,0)(35,11)(45,0)
  \end{picture}
  }

  \put(0,22.5){
  \begin{picture}(60,15)
    \put(0,0){\line(1,0){60}}
    \qbezier(5,0)(20,15)(35,0)
    \qbezier(15,0)(30,15)(45,0)
    \qbezier(25,0)(40,15)(55,0)
  \end{picture}
  }
  
  \put(0,7.5){
  \begin{picture}(60,15)
    \put(0,0){\line(1,0){60}}
    \qbezier(5,0)(15,11)(25,0)
    \qbezier(15,0)(30,15)(45,0)
    \qbezier(35,0)(45,11)(55,0)
  \end{picture}
  }

  \put(75,75){
  \begin{picture}(50,15)
    \put(0,0){\line(1,0){50}}
    \qbezier(5,0)(20,15)(35,0)
    \qbezier(15,0)(25,11)(35,0)
    \qbezier(25,0)(35,11)(45,0)
  \end{picture}
  }

  \put(75,60){
  \begin{picture}(50,15)
    \put(0,0){\line(1,0){50}}
    \qbezier(5,0)(15,11)(25,0)
    \qbezier(15,0)(25,11)(35,0)
    \qbezier(5,0)(25,19)(45,0)
  \end{picture}
  }

  \put(75,45){
  \begin{picture}(50,15)
    \put(0,0){\line(1,0){50}}
    \qbezier(5,0)(25,19)(45,0)
    \qbezier(15,0)(25,11)(35,0)
    \qbezier(25,0)(35,11)(45,0)
  \end{picture}
  }
  
  \put(75,30){
  \begin{picture}(50,15)
    \put(0,0){\line(1,0){50}}
    \qbezier(5,0)(15,11)(25,0)
    \qbezier(5,0)(20,15)(35,0)
    \qbezier(15,0)(30,15)(45,0)
  \end{picture}
  }
  
  \put(75,15){
  \begin{picture}(50,15)
    \put(0,0){\line(1,0){50}}
    \qbezier(5,0)(15,11)(25,0)
    \qbezier(15,0)(25,11)(35,0)
    \qbezier(15,0)(30,15)(45,0)
  \end{picture}
  }
  
  \put(75,0){
  \begin{picture}(50,15)
    \put(0,0){\line(1,0){50}}
    \qbezier(5,0)(20,15)(35,0)
    \qbezier(15,0)(30,15)(45,0)
    \qbezier(25,0)(35,11)(45,0)
  \end{picture}
  }
  
  \put(140,37.5){
  \begin{picture}(40,15)
   \put(0,0){\line(1,0){40}}
   \qbezier(5,0)(20,15)(35,0)
   \qbezier(5,0)(15,11)(25,0)
   \qbezier(15,0)(25,11)(35,0)
  \end{picture}
  }
  
\end{picture}
$$
The dual bracket diagrams are
\begin{center}
\begin{tabular}{ccc}
  \begin{tabular}{c}
    $[x_{1},x_{6}]\cdot [x_{2},x_{4}]\cdot [x_{3},x_{5}]$\\
    $[x_{1},x_{5}]\cdot [x_{2},x_{4}]\cdot [x_{3},x_{6}]$\\
    $[x_{1},x_{4}]\cdot [x_{2},x_{6}]\cdot [x_{3},x_{5}]$\\
    $[x_{1},x_{4}]\cdot [x_{2},x_{5}]\cdot [x_{3},x_{6}]$\\
    $[x_{1},x_{3}]\cdot [x_{2},x_{5}]\cdot [x_{4},x_{6}]$\\
  \end{tabular} &
  \begin{tabular}{c}
    $[x_{1},[x_{2},x_{4}]]\cdot [x_{3},x_{5}]$\\
    $[[x_{1},x_{3}],x_{5}]\cdot [x_{2},x_{4}]$\\
    $[x_{1},[x_{3},x_{5}]]\cdot [x_{2},x_{4}]$\\ 
    $[[x_{1},x_{3}],x_{4}]\cdot [x_{2},x_{5}]$\\ 
    $[x_{1},x_{3}]\cdot [[x_{2},x_{4}],x_{5}]$\\
    $[x_{1},x_{4}]\cdot [x_{2},[x_{3},x_{5}]]$
  \end{tabular}&
  \begin{tabular}{c}
   $[[x_{1},x_{3}],[x_{2},x_{4}]]$
  \end{tabular}
\end{tabular} 
\end{center}

\section{Computer calculations}\label{B}
We give here results of computations of the homology bialgebra of the 4 complexes:  $CT_0D^{even}$,  
$CT_0D^{odd}$,  $CTD^{even}$,  $CTD^{odd}$.  These complexes are graded commutative differential Hopf algebras.

I remind the grading $i$ is the complexity (number of edges) of the corresponding diagrams, the grading $j$ is the
number of points. The numbers $i$ and $j$ correspond to the coordinates $p$ and $q$ of the Vassiliev spectral
sequence by the following formulae:
$$
p=-i, \qquad q=di-j.
$$
Hence $p+q=i(d-1)-j$ is the corresponding cohomology degree of the space of long knots in ${\mathbb R}^d$, $d\ge
3$.

Note that the grading $j$ is allways odd, but the grading $i$ is odd in the case of $d$ even and is even in the
case of $d$ odd.

In the tables below one puts nothing in a cell if there are no diagrams of this bigrading in the corresponding
complex. A question symbol "?" without anything else in a cell means that my computer did not manage to compute
this homology group. A question symbol with some extra information in a cell means that this information is not
sure --- big numbers appeared in the computations.

The groups from the bialgebra of chord diagrams are underlined. Their ranks are known till the complexity 
$i\leq 12$.

The results in tables 5-6 follow from the results of tables 1-4.

\vspace{2mm}

\begin{center}
\begin{tabular}{ccc}
{\bf 1.} Homology of  complex $CT_{0}D^{even}$&&{\bf 2.} Homology of  complex $CT_{0}D^{odd}$\\
  \begin{tabular}{|c|c|c|c|c|c|c|}
       \hline
       & \multicolumn{6}{|c|}{\bf Complexity $i$} \\ \cline{2-7}
    {\bf $j$}& \bf 1& \bf 2& \bf 3& \bf 4& \bf 5& \bf 6 \\ \hline
    \bf 1&&&&&& \\ \hline
    \bf 2&&&&&& \\ \hline
    \bf 3&&&&&& \\  \hline
    \bf 4&& $\Z$& $0$&&& \\ \hline
    \bf 5&&& $\Z$& $0$&& \\ \hline
    \bf 6&&& $\Z$& $\Z_{10}$& $0$& \\ \hline
    \bf 7&&&& $\Z_2^2$& $\Z$& $0$ \\ \hline
    \bf 8&&&& $\Z$& $\Z_2$& ? \\  \hline
    \bf 9&&&&& $\Z^2\oplus \Z_2$?& ? \\ \hline
    \bf 10&&&&& $\Z^3$?& ? \\ \hline
    \bf 11&&&&&& ? \\ \hline
    \bf 12&&&&&& ? \\ \hline
  \end{tabular}
&&
 \begin{tabular}{|c|c|c|c|c|c|c|}
   \hline
   & \multicolumn{6}{|c|}{\bf Complexity $i$}\\ \cline{2-7}
   $j$& \bf 1& \bf 2& \bf 3& \bf 4& \bf 5& \bf 6 \\ \hline
   \bf 1&&&&&& \\ \hline
   \bf 2&&&&&& \\ \hline
   \bf 3&&&&&& \\ \hline
   \bf 4&& \underline{$\Z$}& $0$&&& \\ \hline
   \bf 5&&& $\Z$& $0$&& \\ \hline
   \bf 6&&& \underline{$\Z$}& $\Z_2$& $0$& \\  \hline
   \bf 7&&&& $\Z^2$& $\Z_2$& $0$ \\  \hline
   \bf 8&&&& \underline{\,$\Z^3$}& $0$& ? \\  \hline
   \bf 9&&&&& $\Z^4$& ? \\  \hline
   \bf 10&&&&& \underline{\,$\Z^4$}& ? \\  \hline
   \bf 11&&&&&& ? \\  \hline
   \bf 12&&&&&& \underline{\,$\Z^9$} \\  \hline
 \end{tabular}
\end{tabular}
\end{center}

\vspace{4mm}

\begin{center}
\begin{tabular}{ccc}
 {\bf 3.} Homology of  complex $CTD^{even}$&&  {\bf 4.} Homology of complex $CTD^{odd}$\\
   \begin{tabular}{|c|c|c|c|c|c|}
   \hline
   &\multicolumn{5}{|c|}{\bf Complexity $i$}\\ \cline{2-6}
   $j$& \bf 1& \bf 2& \bf 3& \bf 4& \bf 5 \\ \hline
   \bf 1&&&&&  \\  \hline
   \bf 2& $\Z$&&&& \\ \hline
   \bf 3&& $\Z_2$&&& \\ \hline
   \bf 4&& $\Z$& $\Z_3$&& \\  \hline
   \bf 5&&& $\Z\oplus\Z_2$& $\Z_2$& \\  \hline
   \bf 6&&& $\Z^2$& $\Z_{30}$& $\Z_5$ \\  \hline
   \bf 7&&&& $\Z\oplus\Z_2^4$& ?    \\    \hline
   \bf 8&&&& $\Z^2$& ?   \\  \hline
   \bf 9&&&&& ? \\  \hline
   \bf 10&&&&& ? \\  \hline
 \end{tabular}
&&
 \begin{tabular}{|c|c|c|c|c|c|}
   \hline
   &\multicolumn{5}{|c|}{\bf Complexity $i$} \\  \cline{2-6}
   $j$& \bf 1& \bf 2& \bf 3& \bf 4& \bf 5 \\  \hline
   \bf 1&&&&& \\  \hline
   \bf 2& \underline{$\Z$}&&&& \\ \hline
   \bf 3&& $\Z$&&&  \\ \hline
   \bf 4&& \underline{\,$\Z^2$}& $0$&&  \\  \hline
   \bf 5&&& $\Z^2$& $\Z_2$&  \\  \hline
   \bf 6&&& \underline{\,$\Z^3$}& $\Z_2$& $0$ \\  \hline
   \bf 7&&&& $\Z^5$& ?  \\  \hline
   \bf 8&&&& \underline{\,$\Z^6$}& ?  \\  \hline
   \bf 9&&&&& ?   \\  \hline
   \bf 10&&&&& \underline{\,$\Z^{10}$}  \\  \hline
 \end{tabular}
\end{tabular}
\end{center}

\vspace{4mm}

\begin{center}
\begin{tabular}{ccc}
 {\bf 5.} Primitive generators of the&&{\bf 6.} Primitive generators of the\\
  homology bialgebra of $DHATD_{\Q}^{even}$&& homology bialgebra  of $DHATD_{\Q}^{odd}$\\
   \begin{tabular}{|c|c|c|c|c|c|c|}
  \hline
  & \multicolumn{6}{|c|}{\bf Complexity $i$}\\  \cline{2-7}
  $j$& \bf 1& \bf 2& \bf 3& \bf 4& \bf 5& \bf 6 \\  \hline
  \bf 1&&&&&&  \\  \hline
  \bf 2& 1&&&&&  \\  \hline
  \bf 3&& 0&&&&  \\  \hline
  \bf 4&& 1& 0&&&  \\  \hline
  \bf 5&&& 1& 0&&  \\  \hline
  \bf 6&&& 1& 0& 0&   \\  \hline
  \bf 7&&&& 0& 1& 0   \\   \hline
  \bf 8&&&& 0& 0& ?   \\   \hline
  \bf 9&&&&& 1? & ?   \\   \hline
  \bf 10&&&&& 2?& ?   \\   \hline
  \bf 11&&&&&& ?      \\  \hline
  \bf 12&&&&&& ?    \\  \hline
 \end{tabular}
&&
\begin{tabular}{|c|c|c|c|c|c|c|}
   \hline
   & \multicolumn{6}{|c|}{\bf Complexity $i$} \\  \cline{2-7}
   $j$& \bf 1& \bf 2& \bf 3& \bf 4& \bf 5& \bf 6 \\  \hline
   \bf 1&&&&&&  \\  \hline
   \bf 2& \underline 1&&&&&  \\  \hline
   \bf 3&& 1&&&&  \\  \hline
   \bf 4&& \underline 1& 0&&&  \\  \hline
   \bf 5&&& 1& 0&&  \\  \hline
   \bf 6&&& \underline 1& 0& 0&  \\  \hline
   \bf 7&&&& 2& 0& 0  \\  \hline
   \bf 8&&&& \underline 2& 0& ?  \\  \hline
   \bf 9&&&&& 3& ?   \\   \hline
   \bf 10&&&&& \underline 3& ?  \\   \hline
   \bf 11&&&&&& ?    \\   \hline
   \bf 12&&&&&& \underline 5    \\   \hline
 \end{tabular}
\end{tabular}
\end{center}

To obtain the table of primitive generators of the homology bialgebra of the complex of $T_0$-diagrams (in both
even and odd cases) one should remove from the above tables the content of two cells  $i=1$, $j=2$ and $i=2$, $j=3$. The content of
all the other cells will be the same.

\begin{center}
 \begin{tabular}{cccc}
 {\bf 7.} Homology ranks of &&&{\bf 8.}  Homology ranks of  \\
 complex  $CT_0D_{\Z_{2}}$&&&complex of  $CTD_{\Z_{2}}$\\
   \begin{tabular}{|c|c|c|c|c|c|c|c|}
    \hline
    &\multicolumn{7}{|c|}{\bf Complexity $i$} \\   \cline{2-8}
    $j$& \bf 1& \bf 2& \bf 3& \bf 4& \bf 5& \bf 6& \bf 7  \\  \hline
    \bf 1&&&&&&&  \\   \hline
    \bf 2&&&&&&&  \\   \hline
    \bf 3&&&&&&& \\   \hline
    \bf 4&& \underline{1}& 0&&&& \\  \hline
    \bf 5&&& 1& 0&&& \\  \hline
    \bf 6&&& \underline{1}& 1&0&&  \\  \hline
    \bf 7&&&& 3&1& 0&  \\  \hline
    \bf 8&&&& \underline{3}&1& 1& 0   \\  \hline
    \bf 9&&&&& 4&3& ?   \\   \hline
    \bf 10&&&&& \underline{4}& 7& ?  \\   \hline
    \bf 11&&&&&& 13& ?  \\   \hline
    \bf 12&&&&&&  \underline{9}& ?  \\  \hline
    \bf 13&&&&&&& ? \\   \hline
    \bf 14&&&&&&& \underline{14}\\   \hline
  \end{tabular}
&&&
  \begin{tabular}{|c|c|c|c|c|c|c|}
    \hline
    & \multicolumn{6}{|c|}{\bf Complexity $i$} \\  \cline{2-7}
    $j$& \bf 1& \bf 2& \bf 3& \bf 4& \bf 5& \bf 6 \\  \hline
    \bf 1&&&&&&  \\   \hline
    \bf 2& \underline{1}&&&&& \\ \hline
    \bf 3&& 1&&&& \\ \hline
    \bf 4&& \underline{2}& 0&&&  \\  \hline
    \bf 5&&& 2& 1&&  \\  \hline
    \bf 6&&& \underline{3}& 2& 0&  \\  \hline
    \bf 7&&&& 6& 2& 0  \\  \hline
    \bf 8&&&& \underline{6}& 4& 2  \\  \hline
    \bf 9&&&&& 11& ?   \\  \hline
    \bf 10&&&&& \underline{10}& ? \\  \hline
    \bf 11&&&&&& ? \\    \hline
    \bf 12&&&&&& \underline{19} \\   \hline
    \multicolumn{7}{c}{}\\
    \multicolumn{7}{c}{}\\
  \end{tabular}
\end{tabular}
\end{center}

\vspace{4mm}

\begin{center}
\begin{tabular}{cc}
{\bf 9.} Ranks of complex  $CT_0D$&{\bf 10.} Ranks of complex  $CTD$\\
 \scriptsize
 \begin{tabular}{|c|c|c|c|c|c|c|c|}
   \hline
  & \multicolumn{7}{|c|}{\bf Complexity $i$} \\ \cline{2-8}
  $j$& \bf 1& \bf 2& \bf 3& \bf 4& \bf 5& \bf 6& \bf 7 \\  \hline
  \bf 1&&&&&&& \\ \hline
  \bf 2&&&&&&& \\ \hline
  \bf 3&&&&&&& \\ \hline
  \bf 4&& \underline 1& 1&&&& \\ \hline
  \bf 5&&& 6& 6&&& \\ \hline
  \bf 6&&& \underline{5}& 39& 34&& \\ \hline
  \bf 7&&&& 68& 284& 216&  \\  \hline
  \bf 8&&&& \underline{36}& 771& 2301& 1566 \\  \hline
  \bf 9&&&&& 850& 8634& 20624 \\ \hline
  \bf 10&&&&& \underline{329}& 14835& 100154  \\ \hline
  \bf 11&&&&&& 11940& 237840 \\  \hline
  \bf 12&&&&&& \underline{3655}& 297620 \\ \hline
  \bf 13&&&&&&& 188720 \\   \hline
  \bf 14&&&&&&& \underline{47844} \\   \hline
 \end{tabular}
&
 \scriptsize
 \begin{tabular}{|c|c|c|c|c|c|c|c|}
  \hline
  &\multicolumn{7}{|c|}{\bf Complexity $i$}  \\  \cline{2-8}
  $j$& \bf 1& \bf 2& \bf 3& \bf 4& \bf 5& \bf 6& \bf 7 \\  \hline
  \bf 1&&&&&&& \\  \hline
  \bf 2& \underline 1&&&&&&  \\  \hline
  \bf 3&& 2&&&&&  \\  \hline
  \bf 4&& \underline 3& 6&&&&  \\  \hline
  \bf 5&&& 20& 24&&&   \\  \hline
  \bf 6&&& \underline{15}& 130& 120&&  \\  \hline
  \bf 7&&&& 210& 924& 720&  \\  \hline
  \bf 8&&&& \underline{105}& 2380& 7308& 5040  \\  \hline
  \bf 9&&&&& 2520& 26432& 64224 \\  \hline
  \bf 10&&&&& \underline{945}& 44100& 303660 \\  \hline
  \bf 11&&&&&& 34650& 705320 \\ \hline
  \bf 12&&&&&& \underline{10395}& 866250 \\ \hline
  \bf 13&&&&&&& 540540 \\  \hline
  \bf 14&&&&&&& \underline{135135} \\ \hline
 \end{tabular}
\end{tabular}
\end{center}

\section{$q$-combinatorics}\label{C}
Let us denote by ${{n+m}\choose n}_q$ the following polynomial over $q$:
$$
{{n+m}\choose n}_q=\sum_{s\in S(n,m)}q^{\alpha(s)},\eqno(\numb\label{eqc0})
$$
where $S(n,m)$ is the set of all shuffles of a cardinality $n$ set $\underbrace{x,x,\dots,x}_{\text{$n$ times}}$
with a cardinality $m$ set $\underbrace{y,y,\dots,y}_{\text{$m$ times}}$; $\alpha(s)$ is the minimal number of
transpositions that one needs to obtain shuffle $s$ from the initial shuffle $\underbrace{x,x,\dots,x}_{\text{$n$
times}},\underbrace{y,y,\dots,y}_{\text{$m$ times}}$.

It can be easily verified that
$$
{{n+m+1}\choose n}_q={{n+m}\choose n}_q+q^{m+1}{{n+m}\choose{n-1}}_q;
\eqno(\numb)\label{eqc1}
$$
$$
\sum_{i=0}^{n}{{n}\choose i}_qq^{\frac{i(i-1)}2}x^i=(1+x)(1+qx)\dots(1+q^{n-1}x). \eqno(\numb)\label{eqc2}
$$

Denote by $n_q$ the polynomial $1+q+q^2+\dots+q^{n-1}$, and by $n_q!=1_q\cdot 2_q\cdot\ldots\cdot n_q$.

One has
$$
{{n+m}\choose n}_q=\frac{(n+m)_q!}{n_q!\, m_q!}. \eqno(\numb\label{eqc3})
$$

In this paper we use only two situations $q=\pm 1$. In the case $q=-1$ the last formula can not be applied because of
division by zero. To find numbers ${{n+m}\choose n}_{-1}$ one should use the formula~\eqref{eqc2}. If $n$ is even:
\begin{multline}
\sum_{i=0}^{i=n}{{n}\choose i}_{-1}(-1)^{\frac{i(i-1)}2}x^i=(1+x)(1-x)\dots(1+x)(1-x)=\\
=(1-x^2)^{\frac n2}=\sum_{j=0}^{\frac
n2}(-1)^j{{\frac n2}\choose j}x^{2j}.\label{eqc4}
\end{multline}
If $n$ is odd:
\begin{multline}
\sum_{i=0}^{i=n}{{n}\choose i}_{-1}(-1)^{\frac{i(i-1)}2}x^i=(1+x)(1-x)\dots(1+x)(1-x)(1+x)=\\
=(1-x^2)^{\frac {n-1}2}(1+x)=
\sum_{i=0}^n(-1)^{\left[\frac i2\right]} {{\frac {n-1}2}\choose \left[\frac i2\right]}x^i.\label{eqc5}
\end{multline}

\eqref{eqc4},~\eqref{eqc5} imply:
$$
{{n+m}\choose n}_{-1}=
\begin{cases}
0,&\text{$n$ and $m$ are odd;}\\
\left({\left[\frac{n+m}2\right]}\atop{\left[ \frac n2\right]}\right), &\text{otherwise.}
\end{cases}\eqno(\numb\label{eqc6})
$$

\bigskip

Denote by $\{  n_1,n_2,\dots,n_k\}_q$ the following polynomial over $q$
$$
\{  n_1,n_2,\dots,n_k\}_q=\sum_{s\in S(n_1,n_2,\dots,n_k)}q^{\alpha(s)},\eqno(\numb\label{eqc6'})
$$
where $S(n_1,n_2,\dots,n_k)$ is the set of all shuffles of  sets $X_1=\{\underbrace{x_1,x_1,\dots,x_1}_{\text{$n_1$
times}}\}$,   $\{\underbrace{x_2,x_2,\dots,x_2}_{\text{$n_2$ times}}\}$,$\dots$,
$\{\underbrace{x_k,x_k,\dots,x_k}_{\text{$n_k$ times}}\}$; $\alpha(s)$ is the minimal number of transpositions that
one needs to obtain shuffle $s$ from the initial shuffle $x_1,x_1,\dots,x_1,x_2,\dots,x_2,\dots,x_k$.

This expression is a generalization of~\eqref{eqc0}:
$$
\{  n,m\}_q={ {n+m}\choose n }_q.
$$

One has
$$
\{  n_1,n_2,\dots,n_k\}_q=\frac{(n_1+n_2+\dots+n_k)_q!}{(n_1)_q!\,(n_2)_q!\dots
(n_k)_q!}.\eqno(\numb\label{eqc7})
$$
This expression is symmetric and  satisfies the following identity:
\begin{multline}
\{  n_1,n_2,\dots,n_i,\dots,n_k\}_q\cdot\{  n_{i1},n_{i2},\dots,n_{im}\}_q=\\
\{ 
n_1,\dots,n_{i-1},n_{i1},n_{i2},\dots,n_{im},n_{i+1},\dots,n_k\}_q, \label{eqc8}
\end{multline}
where $n_i=n_{i1}+n_{i2}+\dots+n_{im}$.

The formula~\eqref{eqc7} can not be used in the case $q=-1$. But in this case~\eqref{eqc8} and~\eqref{eqc6} imply:
\begin{equation}
\{  n_1,n_2,\dots,n_k\}_{-1}=
\begin{cases}
0, \, \text{\footnotesize if at least two numbers $n_i$ and $n_j$, $\{i,j\}\subset\{1,\dots,k\}$, are odd;}\\
\{  \left[\frac{n_1}2\right],\left[\frac{n_2}2\right],\dots,\left[\frac{n_k}2\right]\},\qquad\text{\footnotesize otherwise.}
\end{cases}
\label{eqc9}
\end{equation}

\section{Homology bialgebras of $\ZZ^{odd}$, $\ZZ^{even}$}\label{D}
In this section we describe the homology bialgebras of $\ZZ^{odd}$, $\ZZ^{even}$ over $\Z_{p}$, $p$ being any prime,
and over $\Q$. 
These results are well known, see~\cite{F,Vain,V0,Ma}. Calculation of
the Bochschtein homomorphism shows that $\Z$ homology groups of $\ZZ^{odd}$, $\ZZ^{even}$ have no higher torsions
$\Z_{p^k}$ with $k\ge 2$ ($p$ being any prime), see~\cite{F,Vain,V0,Ma}. This permit to find the corresponding homology
groups for any ring $\kk$ of coefficients.
Complexes $\ZZ^{odd}\otimes\kk$, $\ZZ^{even}\otimes\kk$ will be denoted by $\ZZ^{odd}_\kk$,
$\ZZ^{even}_\kk$. Since $\ZZ^{odd}_{\Z_2}\simeq\ZZ^{even}_{\Z_2}$, this bialgebra will be denoted by $\ZZ_{\Z_2}$.
The homology bialgebras of $\ZZ^{odd}_\kk$, $\ZZ^{even}_\kk$, $\ZZ_{\Z_2}$ will be denoted by $\HH^{odd}_\kk$,
$\HH^{even}_\kk$, $\HH_{\Z_2}$.

\smallskip
First, we give necessary definitions. Most part of them are standard.

\begin{definition}\label{dd1}
{\rm A bialgebra is called {\it polynomial} if it is polynomial as algebra and  its 
generators are primitive (if $2\neq0$ in $\kk$, all generators must be of even degree).} \nbox
\end{definition}

Note that over $\Z_p$ the $p$-th power $x^p$ of any primitive element $x$ is always a primitive element: $\Delta
x^p=x^p\otimes 1+1\otimes x^p$.

\begin{definition}\label{dd2}
{\rm A bialgebra is called {\it exterior} if it is exterior as algebra and all its  generators are primitive (if 
$2\neq0$ in $\kk$, all generators must be of odd degree).} \nbox
\end{definition}

\begin{definition}\label{dd3}
{\rm We call  {\it graded polynomial bialgebra} any tensor product of a polynomial bialgebra with an exterior
bialgebra.} \nbox
\end{definition}

Special case is when the ground ring is of characteristic 2. In this situation generators of a polynomial (resp. exterior)
bialgebra can be odd (resp. even) elements. Actually this case makes us emphazise the difference between 
polynomial and exterior  bialgebras.

\begin{definition}\label{dd4}
{\rm A bialgebra $\Gamma$ that is dual to a polynomial bialgebra $A$ is called {\it bialgebra of divided powers}.
The   space in $\Gamma$ dual to the space of generators of $A$ will be called {\it space of divided powers
generators}}. \nbox
\end{definition}

As example consider a bialgebra $\Gamma_\kk(x)$ of divided powers that is dual to a polynomial bialgebra with only
one generator, \cf~\cite{Br}. $\kk$ designates the ground ring. The space of $\Gamma_\kk(x)$ is linearly spanned by the elements
$1=x^{(0)},x=x^{(1)},x^{(2)},x^{(3)},\dots$. If $\kk\neq \Z_2$ all these elements are even. Multiplication and
comultiplication of the elements are given as follows:
\begin{gather}
x^{(k)}\cdot x^{(l)}={ {k+l}\choose k} x^{(k+l)};\label{eqd4'}\\
\Delta(x^{(k)})=\sum_{i=0}^kx^{(i)}\otimes x^{(k-i)}.\label{eqd4''}
\end{gather}
If $\kk=\Z_p$, $p$ being any prime, $\Gamma_\kk(x)$ is a commutative algebra with generators $y_k=x^{(p^k)}$,
$k=0,1,2,\dots$, the only relations are $(y_k)^p=0$. Note that only the generator $y_1=x^{(1)}$ is primitive.

Bialgebra $\Gamma_\Z(x)$ is isomorphic to  a $\Z$-subbialgebra of the polynomial bialgebra $\Q[x]$ $\Z$-spanned by
the elements $x^{(k)}=\frac{x^k}{k!}$. Obviously, $\Gamma_\kk(x)\simeq\Gamma_\Z(x)\otimes\kk$. Note also that
$\Gamma_\Q(x)\simeq \Q[x]$.

\medskip

A bialgebra dual to an exterior bialgebra is also an exterior bialgebra.

\begin{definition}\label{dd5}
{\rm A bialgebra that is dual to a graded polynomial bialgebra is called {\it graded bialgebra of divided
powers}.} \nbox
\end{definition}

A graded bialgebra of divided powers is a tensor product of a bialgebra of divided powers with an exterior
bialgebra.

\vspace{5mm}

For any element $x$ in $\ZZ^{odd}$ or in $\ZZ^{even}$ one can define its {\it divided powers}:
$$
x^{\langle n\rangle }:=\langle \underbrace{x,x,\dots,x}_{\text{$n$ times}}\rangle,\,\, n=0,1,2,\dots.
$$
In the cases

1) $\kk=\Z_2$,

2) $deg(x)$ is even and $\kk$ being an arbitrary ring of coefficients,

\noindent  multiplication and comultilication identities~\eqref{eqd4'},~\eqref{eqd4''} are satisfied by the elements
$x^{\langle n\rangle }$, $n=0,1,2,\dots$ and the following formula holds:
$$
\partial x^{\langle n\rangle }=\partial x\cdot x^{\langle n-1\rangle }.
$$
Due to the last formula, divided powers are homological operations in the described situations. Note, this is
also true for $DHATD$, $DHAT_0D$, $DHAT_*D$. In fact the latter bialgebras are so called {\it divided systems}, 
\cf~\cite[Chapter V, p.~124]{Br}.

\medskip

The element $Z_1\in\ZZ^{odd}$ is even. Denote by $\iota$ the map
$$
\iota:\Gamma_\Z(x)\to\ZZ^{odd},
$$
that sends
$$
\iota:x^{(k)}\mapsto (Z_1)^{\langle k\rangle }.
$$
This map is a morphism of differential Hopf algebras (we define the differential in $\Gamma_\Z(x)$ as zero).

Denote by
$$
I_\Z:\ZZ^{even}\to\ZZ^{odd},
$$
the map that sends
$$
I_\Z:\langle Z_{k_1},Z_{k_2},\dots,Z_{k_\ell}\rangle \mapsto\langle Z_{2k_1},Z_{2k_2},\dots,Z_{2k_\ell}\rangle .
$$
$I_\Z$ is a morphism of differential Hopf algebras. To see this one should compare  formulas~\eqref{eq114}
and~\eqref{eq115}.

We will formulate several assertions in order to describe the homology bialgebras of $\ZZ^{even}$, $\ZZ^{odd}$.
The most part of these lemmas are wether well known or are reformulations of well known results.

\begin{lemma}\label{ld5}
The composition map
$$
\xymatrix{
\Gamma_\Z(x)\otimes\ZZ^{even}\ar[r]^-{\iota\otimes I_\Z}&\ZZ^{odd}\otimes\ZZ^{odd}\ar[r]^-\mu&\ZZ^{odd} 
} \eqno(\numb)\label{eqd6}
$$
is a quasi-isomorphism of differential Hopf algebras. \nbox
\end{lemma}

\begin{corollary}\label{cd5'}
Since~\eqref{eqd6} is a quasi-isomorphism over $\Z$, it stays quasi-isomorphism for any ring $\kk$ of
coefficients. In the special case $\kk=\Z_2$ quasi-isomorphism~\eqref{eqd6} implies the quasi-isomorphism:
$$
\xymatrix{
\Gamma_{\Z_2}(x)\otimes\ZZ_{\Z_2}\ar[rr]^-{\mu\circ(\iota_{\Z_2}\otimes I_{\Z_2})}&&\ZZ_{\Z_{2}},
}\eqno(\numb)\label{eqd7}
$$
where $\iota_{\Z_2}=\iota\otimes\Z_2$, $I_{\Z_2}=I_Z\otimes\Z_2$.
\end{corollary}

Applying this quasi-isomorphism infinitely many times we obtain:

\begin{theorem}\label{td7}
Differential Hopf algebra $\ZZ_{\Z_2}$ is formal. The homology bialgebra $\HH_{Z_2}$ of $\ZZ_{\Z_2}$ is a
bialgebra of divided powers. Divided powers generators of $\HH_{Z_2}$ are elements $Z_{2^k}$, $k=0,1,2,\dots$. As
algebra $\HH_{Z_2}$ is generated by $x_{k;n}=(Z_{2^k})^{\langle 2^n\rangle }$, $k,n=0,1,2,\dots$. The only relations are
$(y_{k;n})^2=0$ for all possible $k$ and $n$. \nbox
\end{theorem}

 A differential Hopf algebra is said {\it formal} if it is quasi-isomorphic to its differential Hopf
subalgebra, such that the restriction of the differential to this subalgebra is zero.

\medskip

For any ring $\kk$ of coefficients Lemma~\ref{ld5} implies that
$\HH^{odd}_\kk\simeq\HH^{even}_{\kk}\otimes\Gamma_\kk(x)$. So, let us concentrate on complex $\ZZ^{even}$.

Suppose $\kk=\Z_{p}$, $p$ being any prime.
Denote by $\ZZ^{even}_{<p}$ a subspace of $\ZZ^{even}_{\Z_p}$ linearly spanned by all the elements
$\langle Z_{k_1},Z_{k_2},\dots,Z_{k_l}\rangle $, such that $k_i<p$, $i=1,2,\dots,l$. It can be easily
verified that $\ZZ^{even}_{<p}$ is a differential Hopf subalgebra of $\ZZ^{even}_{\Z_p}$. Really, multiplication
and comultiplication preserve this space. Since ${{a+b}\choose a}\equiv 0\mod p$ if $a,b<p$ and $a+b\ge p$, the
differential does also preserve $\ZZ^{even}_{<p}$. By abuse of the langauge the inclusion
$\ZZ^{even}_{<p}\hookrightarrow\ZZ^{even}_{\Z_p}$ will be denoted by $\iota$.

Denote by
$$
I_{\Z_p}:\ZZ^{even}_{\Z_p}\to\ZZ^{even}_{\Z_p}
$$
the map that sends
$$
I_{\Z_p}:\langle Z_{k_1},Z_{k_2},\dots,Z_{k_\ell}\rangle \mapsto\langle Z_{pk_1},Z_{pk_2},\dots,Z_{pk_\ell}\rangle .
$$
$I_{\Z_p}$ is a morphism of differential Hopf algebras. To see this, note that ${{pa+pb}\choose pa}\equiv
{{a+b}\choose a}\mod p$ for any $a,b\in\N$.

\begin{lemma}\label{ld8}
For any prime $p$, the composition map
$$
\xymatrix{
\ZZ^{even}_{<p}\otimes\ZZ^{even}_{\Z_p}\ar[r]^-{\iota\otimes I_{\Z_p}}&\ZZ^{even}_{\Z_p}\otimes\ZZ^{even}_{\Z_p}
\ar[r]^-\mu&\ZZ^{even}_{\Z_p}
}\eqno(\numb)\label{eqd9}
$$
is a quasi-isomorphism of differential Hopf algebras. \nbox
\end{lemma}

$\mu$ in~\eqref{eqd9} designates as usual the shuffle multiplication.

In the case $p=2$ Lemma~\ref{ld8} is equivalent to Corollary~\ref{cd5'}.

From now on let $p>2$.

\begin{lemma}\label{ld10}
Any cycle in $\ZZ^{even}_{<p}$, $p$ being any odd prime, is homologous up to a coefficient to one of the cycles:
\begin{gather}
\langle \underbrace{Z_1,Z_{p-1},\dots,Z_1,Z_{p-1}}_{\text{$n$ {\rm times}
$Z_1,Z_{p-1}$}}\rangle ,\,\, n\ge 0,\label{eqd11}\\
\langle \underbrace{Z_1,Z_{p-1},\dots,Z_1,Z_{p-1}}_{\text{$n$ {\rm times}
$Z_1,Z_{p-1}$}},Z_1\rangle ,\,\, n\ge 0. \label{eqd12}
\end{gather}
$\Box$
\end{lemma}

\begin{remark}\label{rd13}
{\rm
\begin{gather*}
\langle \underbrace{Z_1,Z_{p-1},\dots,Z_1,Z_{p-1}}_{\text{$n$ {\rm times}
$Z_1,Z_{p-1}$}}\rangle =\langle Z_1,Z_{p-1}\rangle ^{\langle n\rangle };\\
\langle \underbrace{Z_1,Z_{p-1},\dots,Z_1,Z_{p-1}}_{\text{$n$ {\rm times}
$Z_1,Z_{p-1}$}},Z_1\rangle =\langle Z_1,Z_{p-1}\rangle^{\langle n\rangle }*Z_1,
\end{gather*}
$*$ designates the shuffle product.\nbox
}
\end{remark}

Applying Lemmas~\ref{ld8},~\ref{ld10} and Remark~\ref{rd13} infinitely many times and also applying Lemma~\ref{ld5} we obtain:

\begin{theorem}\label{td14}
1) The homology bialgebra $\HH^{even}_{\Z_p}$ of $\ZZ^{even}_{\Z_p}$, $p$ being any odd prime, is a graded
bialgebra of divided powers. Its divided powers generators are $y_k=\langle Z_{p^k},Z_{(p-1)p^k}\rangle $,
$k=0,1,2,\dots$, its exterior generators are $z_k=Z_{p^k}$, $k=0,1,2,\dots$. As algebra $\HH^{even}_{\Z_p}$
is a graded commutative algebra generated by even elements $y_{k;n}=(y_k)^{\langle p^n\rangle }$, $k,n=0,1,2,\dots$ and odd
elements $z_k$, $k=0,1,2,\dots$. The only relations are $(y_{k;n})^p=0$ for all possible $k$ and $n$.

2) The homology bialgebra $\HH^{odd}_{\Z_p}$ of $\ZZ^{odd}_{\Z_p}$, $p$ being any odd prime, is a graded bialgebra
of divided powers. Its divided powers generators are $x=Z_1$ and $y_k=\langle Z_{2p^k},Z_{2(p-1)p^k}\rangle $,
$k=0,1,2,\dots$, its exterior generators are $z_k=Z_{2p^k}$, $k=0,1,2,\dots$. As algebra $\HH^{odd}_{\Z_p}$
is a graded commutative algebra generated by even elements $x_k=x^{\langle p^k\rangle }$, $k=0,1,2,\dots$,
$y_{k;n}=(y_k)^{\langle p^n\rangle }$, $k,n=0,1,2,\dots$ and odd elements $z_k$, $k=0,1,2,\dots$. The only relations are $(x_k)^p=0$, $(y_{k;n})^p=0$
for all possible $k$ and $n$. \nbox
\end{theorem}

Finally we formulate a theorem describing the situation $\kk=\Q$:

\begin{theorem}\label{td15}
1) The homology bialgebra $\HH^{even}_\Q$ of $\ZZ^{even}_\Q$ is an exterior bialgebra, its only generator is
$Z_1$.

2) The homology bialgebra $\HH^{odd}_\Q$ of $\ZZ^{odd}_\Q$ is a graded polynomial bialgebra with one even
generator $Z_1$ and one odd generator $Z_2$. \nbox
\end{theorem}

\section{Dyer-Lashof operations and explicit formulas for inclusion of $H_{*}(\Omega^2S^{d-1})$ in Hochschild
homology  of $\Poiss$, $\Gerst$}\label{E}

In previous section we described how the cohomology bialgebra of $\Omega^2S^{d-1}$ is contained in 
the first term of Vassiliev cohomological long knots spectral sequence. The aim of this section is to explicitly describe 
the dual situation: how $H_{*}(\Omega^2S^{d-1})$ is included in the homological Vassiliev spectral sequence, \ie in the 
homology of the Hochschild complexes
$(\Poiss,\partial)$, $(\Gerst,\partial)$.

\medskip

It is already well known that Hochschild cochain complex of an associative algebra can be endowed with an action
of an operad quasi-isomorphic to the chain operad of small squares. This statement is named 
"Deligne's conjecture". Over $\Z$ this result is due to J.E.~McClure and J.H.~Smith, \cf~\cite{McSm}, in characteristic 
zero there are several proofs, \cf~\cite{K2,Tam1,Tam2,Vor}. Actually the  operad considered 
by J.E.~McClure and J.H.~Smith consists of brace operations, see Definition~\ref{d811}, and
of an associative mutiplication that satisfy some composition and differential relations. It means that their proof works 
as well in the case of Hochschild complex $(\O,\partial)$ of an operad $\O$, see Section~\ref{ss83}.
This result implies that in Hochschild homology one can define the same homological operations as for 
double loop spaces.  Homological operations for the iterated loop spaces are well known, \cf~\cite{DyLash,Co}:
In the case of double loops, one has Pontriagin multiplication, Browder operator --- degree one bracket  (we will designate it $[\, .\, ,\, .\, ]$), and also two
non-trivial Dyer-Lashof operations (following F.Cohen we designate them $\xi_{1}$ and $\zeta_{1}$):

\smallskip

Over $\Z_{2}$:
 $$
 \xi_{1}:H_{k}(\Omega^2X,\Z_{2})\to H_{2k+1}(\Omega^2X,\Z_{2}).
 $$
 
 \smallskip
 
 Over $\Z_{p}$, $p$ being any odd prime:
 $$
 \xi_{1}:H_{2k-1}(\Omega^2X,\Z_{p})\to H_{2pk-1}(\Omega^2X,\Z_{p}),
 $$
 $$
 \zeta_{1}:H_{2k-1}(\Omega^2X,\Z_{p})\to H_{2pk}(\Omega^2X,\Z_{p}).
 $$
 
 \medskip
 
$\xi_{1}$ and $\zeta_{1}$ are related via Bochstein homomorphism $\beta$:
$$
\zeta_{1}\varphi=\beta\xi_{1}x- (\ad^{p-1}\varphi)(\beta\varphi).
\eqno(\numb\label{eqe1})
$$
Operator $\ad\varphi$ is the ajoint action $[\varphi,\, .\, ]$.

\medskip

Let $\varphi\in\O$ be a cycle of a Hochschild complex $(\O,\partial)$. Suppose the ground ring $\kk=\Z_{2}$ or $deg(\varphi)$ is odd.
By abuse of the language define
$$
\xi_{1}(\varphi):=\varphi^{[p]}.
\eqno(\numb\label{eqe2})
$$
Suppose $\kk=\Z_p$. Define
$$
\zeta_{1}(\varphi):=-\sum_{i=1}^{p-1}\frac{(p-1)!}{i!(p-i)!}\varphi^{[i]}*\varphi^{[p-i]}\equiv
\sum_{i=1}^{p-1}\frac{(-1)^{i}}{i}\varphi^{[i]}*\varphi^{[p-i]}\mod p.
\eqno(\numb\label{eqe3})
$$
Modulo Remark~\ref{hom}, Proposition~\ref{p92} implies that operations~\eqref{eqe2}, \eqref{eqe3} are homological operations.

\begin{proposition}\label{pe1}
Operations $\xi_{1}$, $\zeta_{1}$ satisfy property~\eqref{eqe1}. \nbox
\end{proposition}

In this proposition $\beta$ is the Bochstein homomorphism, $\ad\varphi=[\varphi,\, .\, ]$ is the adjoint Gerstenhaber bracket
action of $\varphi$.

\medskip

\noindent {\bf Proof of Proposition~\ref{pe1}:} Essentially one needs to prove the following identity:
$$
[\underbrace{\varphi,[\varphi,\ldots,[\varphi}_{p-1},\beta\varphi]\ldots]]\equiv\sum_{i=1}^p
(\ldots((\ldots(\varphi\circ\varphi)\circ\ldots)\beta\varphi)\circ\dots)\circ\varphi\mod p.
\eqno(\numb\label{eqe4})
$$
The $i$-th summand of the right hand-side of~\eqref{eqe4} is obtained by substitution of the $i$-th $\varphi$ by $\beta\varphi$
in $(\dots((\underbrace{\varphi\circ\varphi)\circ\varphi)\dots)\circ\varphi}_{\text{$p$ times}}$.

\medskip

Identity~\eqref{eqe4} follows from the following combinatorial fact:

\begin{lemma}\label{le2}
Let $T$ be any tree with $p$ vertices, $a$ and $b$ be its two vertices, then
$$
O(T_{a})\equiv (-1)^{d(a,b)}O(T_{b})\mod p,$$
where $T_{a}$ (resp. $T_{b}$) is  tree $T$ with a choosen root $a$ (resp. $b$), $d(a,b)$ is the number 
of edges of the only path joining vertex $a$ to vertex $b$. \nbox
\end{lemma}

Proposition~\ref{pe1} confirms the following result.

\begin{theorem}\label{cone3}
Homology operations $\xi_{1}$ and $\zeta_{1}$ are the Dyer-Lashof operations induced
by  the little disc chain action on the Hochschid complex. \nbox
\end{theorem}

To prove this theorem one needs to analyse in detail the McClure-Smith construction, that
is beyond the capabilities of this paper. The proof of Theorem~\ref{cone3} will be given by the author elsewhere. Note, this
problem was already solved for $p=2$ in~\cite{Craig}.

Operations $\xi_{1}$, $\zeta_{1}$ and also the bracket permit to give explicit formulas for the inclusion
of $\Omega^2S^{d-1}$ homology classes as homology classes of Hochschild complexes
$(\Poiss,\partial)$ (case of odd $d$) and $(\Gerst,\partial)$ (case of even $d$).

Let $\varphi\in H_{d-3}(\Omega^2 S^{d-1},\kk)$ denote the image of a generator of $\pi_{d-3}(\Omega^2 S^{d-1})\simeq\Z$
via Hurewicz homomorphism. The following theorems describe the homology bialgebra of $\Omega^2 S^{d-1}$, 
\cf~\cite{AK,Co,CLM,DyLash,May,Milg, Segal,Snaith}, 
and
as consequence of Theorems~\ref{t72}, \ref{t116}, \ref{td7}, \ref{td14}, \ref{td15}, 
Proposition~\ref{p92} and Lemmas~\ref{l105'}, \ref{l1010} provide the desired inclusion formulas (we set $\varphi=[x_{1},x_{2}]$).

\begin{theorem}\label{te3}
For any even $d\geq 4$, bialgebra $H_{*}(\Omega^2S^{d-1},\Q)$ is an exterior bialgebra with the only generator $\varphi$.

For any odd $d\geq 5$ bialgebra $H_{*}(\Omega^2S^{d-1},\Q)$ is a graded polynomial bialgebra with two generators: even generator $\varphi$ and 
odd generator ${\frac 12}[\varphi,\varphi]$. \nbox
\end{theorem}

\begin{theorem}\label{te4}
For any $d\geq 4$, bialgebra $H_{*}(\Omega^2S^{d-1},\Z_{2})$ is a graded polynomial bialgebra with generators $\xi_{1}^k\varphi\in 
H_{2^k(d-2)-1}(\Omega^2S^{d-1},\Z_{2})$, $k=0,1,2,\dots$. \nbox
\end{theorem}

\begin{theorem}\label{te5}
For any even $d\geq 4$, bialgebra $H_{*}(\Omega^2S^{d-1},\Z_{p})$, $p$ being any odd prime, is a graded polynomial bialgebra with exterior generators
$\xi_{1}^k\varphi\in 
H_{p^k(d-2)-1}(\Omega^2S^{d-1},\Z_{p})$, $k=0,1,2,\dots$, and polynomial generators
$\zeta_{1}\xi_{1}^k\varphi=\beta\xi_{1}^{k+1}\varphi\in H_{p^{k+1}(d-2)-2}(\Omega^2S^{d-1},\Z_{p})$, $k=0,1,2,\dots$. \nbox
\end{theorem}

\begin{theorem}\label{te6}
For any odd $d\geq 5$, bialgebra $H_{*}(\Omega^2S^{d-1},\Z_{p})$, $p$ being any odd prime, is a graded polynomial bialgebra with exterior generators
$\xi_{1}^k\psi\in 
H_{2p^k(d-2)-1}(\Omega^2S^{d-1},\Z_{p})$, $k=0,1,2,\dots$, where $\psi={\frac 12}[\varphi,\varphi]$, and polynomial generators
$\varphi\in H_{d-3}(\Omega^2S^{d-1},\Z_{p})$ and  
$\zeta_{1}\xi_{1}^k\psi=\beta\xi_{1}^{k+1}\psi\in H_{2p^{k+1}(d-2)-2}(\Omega^2S^{d-1},\Z_{p})$, $k=0,1,2,\dots$. \nbox
\end{theorem}

\noindent Independent University of Moscow

\noindent Universit\'e Catholique de Louvain (Louvain-la-Neuve)

\noindent {\tt e-mail: turchin@mccme.ru, turchin@math.ucl.ac.be}

\end{document}